\documentclass[12pt]{article}
\usepackage{graphicx}
\graphicspath{{images/}}

\usepackage{amscd}
\usepackage{latexsym}
\usepackage{amsthm}
\usepackage{verbatim}
\usepackage[english]{babel}
\usepackage{amsfonts}
\usepackage[a4paper,margin=0.9in]{geometry}
\usepackage{mathtools}
\usepackage{bbm}
\usepackage{lscape}
\usepackage[T1]{fontenc} 
\usepackage{amsmath,amsbsy,amssymb,latexsym, graphics, graphicx, tabularx}
\usepackage{adjustbox}
\usepackage{lscape}
\usepackage{tikz}
\usepackage{color}
\usepackage{url}
\usepackage{rotating}
\usepackage{caption}
\usepackage{authblk}
\usepackage{enumerate}
\usepackage{soul, xcolor}
\usepackage{multirow}
\usepackage{subfigure}
\usepackage{graphicx}
\usepackage{csquotes}
\usepackage[font=small]{caption}
\usepackage{longtable}
\usepackage[shortlabels]{enumitem}
\usepackage{authblk}
\usepackage{booktabs}

\theoremstyle{plain} 
\newtheorem{Lemma}{Lemma}[section] 
\newtheorem{Theorem}{Theorem}[section] 
\newtheorem{Corollary}{Corollary}[section] 

\theoremstyle{plain} 

\theoremstyle{remark} 
\newtheorem{Remark}{Remark}[section] 

\makeatletter
\def\th@remark{%
	\thm@headfont{\bfseries}%
	\normalfont 
	\thm@preskip\topsep \divide\thm@preskip\tw@
	\thm@postskip\thm@preskip
}
\makeatother

\makeatletter

\@addtoreset{equation}{section}
\makeatother

\usepackage{xr}
\usepackage{hyperref}

\hypersetup{
	colorlinks   = true, 
	urlcolor     = blue, 
	linkcolor    = blue, 
	citecolor   = blue 
}

\sethlcolor{yellow}
\setstcolor{red}

\newcommand{\R}{\mathbb{R}}

\newcommand{\Iff}{\Leftrightarrow}

\newcommand{\Var}{\mbox{Var}}

\newcommand{\cpl}{\mathcal{C}}
\newcommand{\eal}{\mbox{EAL}}
\newcommand{\ert}{\epsilon_T}
\newcommand{\erc}{\epsilon_C}
\newcommand{\datan}{\mathcal{D}_n}

\newcommand{\enq}[1]{\enquote*{#1}}

\title{\bf Quantile regression under dependent censoring with unknown association}
\author[1,2]{Myrthe D'Haen}
\author[2]{Ingrid Van Keilegom}
\author[1]{Anneleen Verhasselt}
\affil[1]{Centre for Statistics, Data Science Institute, Hasselt University, Belgium}
\affil[2]{Research Centre for Operations Research and Statistics, KU Leuven, Belgium}
\date{\today}

\begin{document}
	
	\maketitle
	
	\begin{abstract}
		\noindent 	
		The study of survival data often requires taking proper care of the censoring mechanism that prohibits complete observation of the data. Under right censoring, only the first occurring event is observed: either the event of interest, or a competing event like withdrawal of a subject from the study. The corresponding identifiability difficulties led many authors to imposing (conditional) independence or a fully known dependence between survival and censoring times, both of which are not always realistic. However, recent results in survival literature showed that parametric copula models allow identification of all model parameters, including the association parameter, under appropriately chosen marginal distributions. The present paper is the first one to apply such models in a quantile regression context, hence benefiting from its well-known advantages in terms of e.g. robustness and richer inference results. The parametric copula is supplemented with a likewise parametric, yet flexible, \emph{enriched asymmetric Laplace} distribution for the survival times conditional on the covariates. Its asymmetric Laplace basis provides its close connection to quantiles, while the extension with Laguerre orthogonal polynomials ensures sufficient flexibility for increasing polynomial degrees. The distributional flavour of the quantile regression presented, comes with advantages of both theoretical and computational nature. All model parameters are proven to be identifiable, consistent, and asymptotically normal. Finally, performance of the model and of the proposed estimation procedure is assessed through extensive simulation studies as well as an application on liver transplant data.
	\end{abstract}
	
	\bigskip
	
	\noindent \textbf{Key words:} Quantile regression, Dependent censoring, Copulas, Laguerre polynomials, Survival analysis
	\setcounter{footnote}{0}
	
	\bigskip
	
	\noindent \textbf{Funding:} M. D'Haen is funded by a BOF PhD fellowship at Hasselt University (no. R-12215). I. Van Keilegom gratefully acknowledges funding from the FWO and F.R.S.-FNRS (Excellence of Science programme, project ASTeRISK, grant no. 40007517), and from the FWO (senior research projects fundamental research, grant no. G047524N). A. Verhasselt receives funding from Hasselt University BOF grant no. R-10786 and FWO (junior research projects fundamental research, grant no. G011022N).
	
\section{Introduction} \label{sec:intro}
Studying lifetime distributions is often complicated by the presence of censoring, causing the observation of another event rather than the survival time for some individuals or objects in the data sample. Consequently, a lot of literature has been devoted to survival analysis under the presence of censoring, ranging over a variety of censoring mechanisms. One such mechanism that is frequently encountered is right censoring, where only the minimum of the survival time $T$ and some censoring time $C$ is observed; increasingly complex settings have been considered, from fixed to random right censoring. The latter was initially (and still often is) assumed to be independent of the survival time, possibly after conditioning on a set of covariates $X$. In many contexts, however, such an independence assumption may not be realistic. For instance in medical studies on the time to death of a specific disease, patients may withdraw from the study due to reasons related to their health condition, or die from other diseases with potentially related causes. Such contemplations led to the emergence of survival literature on dependent censoring, where the parametric copula approach is a popular way to overcome the nonidentifiability in a fully nonparametric setting revealed by \cite{tsiatis_nonidentifiability_1975}. Typically, to allow nonparametric estimation of the distribution of the survival time $T$ and censoring time $C$, the copula function capturing their association is assumed to be completely known. This includes both the copula family \emph{and} the parameter measuring the strength of this association, which is often unknown in practice. Moreover, it is a consistent finding in the literature that such models are rather robust to misspecification of the copula family, but sensitive to incorrectly specified association levels \cite{zheng_estimates_1995, huang_regression_2008, chen_semiparametric_2010}. With their recent article, \cite{czado_dependent_2023} opened the door to a new approach: the association parameter within a copula family can be identified from the data at the cost of some assumptions on the margins $T$ and $C$. While both were still parametric in that article, the follow-up paper by \cite{deresa_copula_2024} took this one step further by considering a semiparametric Cox proportional hazards model for $T|X$, conditional on the covariates. The model proposed in the present article is similar in spirit -- preserving the parametric copula with identifiable association parameter, a simple parametric model for $C|X$, and a more flexible model for $T|X$ -- yet does not require the proportional hazards assumption. Rather, the general idea consists in taking a model for $T$ that is still parametric, but can approximate any relevant distribution sufficiently well when the number of parameters is increased. Concretely, this is done by augmenting the degree of the Laguerre polynomials that constitute an important part of the proposed model, the introduction of which is an idea taken from \cite{kreiss_VK_submitted} in the context of quantile regression under (conditionally) independent censoring.

The focus on quantile regression in the latter article is also adopted in our work (at least to some extent, cf. Remark \ref{rem:model}\ref{rem:model:distr.reg} below), and constitutes another important difference with the mentioned work of \cite{czado_dependent_2023} and \cite{deresa_copula_2024}. Quantile regression, introduced in the seminal paper by \cite{koenker_regression_1978}, has several advantages over classical mean regression, see e.g. the book by \cite{koenker_quantile_2005}. For instance, it allows one to do inference for the whole distribution rather than only the mean or median. Furthermore, it is more robust against outliers and enjoys some nice theoretical properties, e.g. its \emph{equivariance to monotone transformations} (Chapter 2 in \cite{koenker_quantile_2005}). This means that for any nondecreasing function $h(\cdot)$ and random variable $Z$, $Q_{h(Z)}(p) = h(Q_{Z}(p))$ for any quantile level $p \in (0,1)$, which generally doesn't hold true for the mean. In other words, the variable of interest can equivalently be studied by considering any transformation that is possibly easier to study, e.g. the unrestricted $\log(T)$ instead of the survival time $T$ itself, restricted to $[0, \infty)$.

\cite{portnoy_censored_2003} moreover provides a motivation, inspired by \cite{koenker_reappraising_2001}, on why quantile regression is interesting especially in a survival context: quantiles can deal with heteroscedasticity and inhomogeneity that are often encountered, and they moreover do so in a natural way, studying the directly interpretable survival times rather than, for instance, the hazard in the popular proportional hazards model. That model additionally has the obvious disadvantage of excluding cases where the proportionality cannot hold, for example when some covariate has a negative effect on part of the population but a positive effect elsewhere; quantile regression is proposed as a natural alternative -- or rather complement -- and potential improvement of the traditional Cox regression \cite{portnoy_censored_2003}. \cite{peng_survival_2008} supplement this reasoning by advocating against the use of classical accelerated failure time models in some cases, in favour of more flexible quantile regression. Quantiles have indeed become a popular regression tool in survival analysis; see e.g. the review paper by \cite{peng_quantile_2021}. Yet, it is much less well-studied under dependent censoring, especially with the flexible association as in the present article.
\bigbreak
The structure of the paper is as follows. Section \ref{sec:lit} discusses some related literature; Section \ref{sec:model} introduces the model and motivates its components. In Section \ref{sec:idf} we discuss identifiability. Asymptotic properties of the estimator are considered in Section \ref{sec:asy}, whereas the numerical estimation procedure is clarified in Section \ref{sec:parest}. The finite sample behaviour of the model is studied in Section \ref{sec:sim} through simulation studies and in Section \ref{sec:realdata} by means of a real data application. Section \ref{sec:concl}, finally, contains some concluding remarks. The online Supplementary material contains mathematical details concerning identifiability and asymptotics, as well as a more in-depth algorithm description and additional simulation scenarios.

\section{Related literature} \label{sec:lit}
Being on the intersection of two important statistical domains, namely (copula models for) dependent censoring and quantile regression, the literature related to at least one of those aspects is extensive; work combining both, on the other hand, is scarce. A by no means exhaustive overview of some related literature on (survival) data subject to censoring is presented here. Section \ref{sec:lit:cop}, with a focus on general copula models, is kept to a minimum; given the quantile orientation in our article, Section \ref{sec:lit:qu}, with a focus on quantiles, dives somewhat deeper into the literature.
\bigbreak
To make the copula literature of Section \ref{sec:lit:cop} slightly more concrete, we already introduce the notion of an (Archimedean) copula, that will be important also in the remainder of our paper. A copula $\cpl(\cdot, \cdot)$ is a bivariate distribution function on the unit square with uniform margins; it is called Archimedean if there exists a generator function $\psi(\cdot)$, continuous, convex and strictly increasing on $[0,1]$, with $\psi(1) = 0$, such that $\cpl(u,v) = \psi^{[-1]}\left(\psi(u) + \psi(v)\right)$ for the pseudo-inverse $\psi^{[-1]}$. (See \cite{nelsen_introduction_2006}, for instance, for an introduction to copulas).

\subsection{Copula models for dependent censoring}\label{sec:lit:cop}
Tsiatis' nonidentifiability result \cite{tsiatis_nonidentifiability_1975} led many authors to considering copula models in order to avoid an independence assumption between $T$ and $C$ that is both untestable and in many cases unrealistic. They rely on Sklar's well-known result that the joint distribution of $(T,C)$ can be uniquely written as the composition of a copula function $\cpl(\cdot, \cdot)$ with the marginal distributions for $T$ and $C$, whenever the latter are continuous \cite{sklar_fonctions_1959}. Early references assumed these copulas to be completely known, including the association strength. \cite{zheng_estimates_1995} initiated this approach by introducing the copula-graphic estimator in a setting without covariates; it was further studied in the case of Archimedean copulas by \cite{rivest_martingale_2001}. \cite{braekers_copula-graphic_2005} then extended this to a fixed design regression case, focussing on the asymptotic properties.

Anticipating on the quantile literature in Section \ref{sec:lit:qu}, we already note that \cite{veraverbeke_regression_2006}, as a follow-up article, inverted the estimator for the distribution function proposed in \cite{braekers_copula-graphic_2005} and studied the properties of the corresponding quantile estimator. Our estimator, too, is based on inverting a distribution rather than directly estimating the quantiles, though specifically developed with quantile estimation in mind (as will also become clear when introducing the model components in Section \ref{sec:model}). Another important difference is that the latter article still belongs to the copula-graphic estimator research line, hence allowing the margins to be modelled nonparametrically. By contrast, we pay the price of parametric margins in order to relax the fixed copula assumption and allow estimation of the association parameter.

Other important references on dependent censoring using copulas with a fixed association are, amongst others, the articles by \cite{huang_regression_2008} and \cite{chen_semiparametric_2010}, as well as those collected in the book by \cite{emura_analysis_2018}. A more comprehensive overview of related copula literature with fixed dependence parameter may be found in \cite{czado_dependent_2023}, that, finally, took a new direction by shifting some of the model flexibility of the margins to that of the copula, as explained in Section \ref{sec:intro}. Subsequently, \cite{deresa_copula-based_2022} added covariates to this approach, as well as left truncation and administrative censoring. The fully parametric assumption on both margins, on the other hand, was replaced by a semiparametric Cox proportional hazards assumption for $T|X$ in \cite{deresa_copula_2024}.

\subsection{Quantile regression under censoring}\label{sec:lit:qu}
There is a vast literature on quantile regression for censored data, but the majority assumes that the response of interest and the censoring are independent either conditional on a set of covariates or even independent at all (without conditioning); see e.g. \cite{wang_locally_2009} for early references of both types of independence, and the more recent review in \cite{peng_quantile_2021}. The first such examples are \cite{powell_least_1984} and \cite{powell_censored_1986}, where censoring is assumed to be always observed and moreover constant. Other well-known articles in this context are \cite{portnoy_censored_2003} and \cite{peng_survival_2008}, both working under the strong global linearity assumption that all quantiles below the one of interest be linear, too. This restriction was relaxed by \cite{wang_locally_2009}, assuming linearity only on the specified quantile level. \cite{wey_censored_2014} subsequently proposed a similar estimator that is moreover able to deal with covariates that are discrete or of (moderately) high dimension. (In our work, the linearity assumption can sometimes be even further relaxed to the existence of \emph{some} linear quantile, not necessarily the one under consideration, cf. Section \ref{sec:model:multiqu}.) More recent references on quantile regression under independent censoring include \cite{leng_quantile_2013}, \cite{de_backer_adapted_2019}, \cite{bravo_semiparametric_2020}, and \cite{ewnetu_two-piece_2023}, amongst others. The latter is closest to our setup, sharing the distributional flavour (Remark \ref{rem:model}\ref{rem:model:distr.reg}) and moreover using two-piece, quantile-oriented parametric distributions as well. This also goes for the work by \cite{gijbels_quantile-based_2019,gijbels_semiparametric_2021}, but in a context without censoring.

Finally, the work of \cite{bottai_laplace_2010} deserves specific mentioning in the context of the present article: the asymmetric Laplace (AL) distribution and its close connection to quantiles is key to both works. The proposal of \cite{bottai_laplace_2010} was, however, refuted in a letter to the editor by \cite{koenker_letter_2011}, exposing the inconsistency of their estimator. Yet, the AL-based research line was continued in recent work by \cite{kreiss_VK_submitted}, extending the AL distributions to a sufficiently flexible family. Their work still belongs to the branch on independent censoring literature, but we transplant the key idea to our dependent context (see also Section \ref{sec:model} below).
\bigbreak
Some steps in the direction of dependent censoring were taken in the work of \cite{peng_nonparametric_2007, peng_competing_2009}, both studying quantiles in a competing risks setting, where the most recent one also includes covariates. They consider nonparametric quantile inference for multiple, mutually censoring events, but they work on the level of the cause-specific cumulative incidence functions rather than those of the latent variables. Whether their so-called \emph{crude} quantities are always preferred over the \emph{net} quantities, in the terminology of \cite{tsiatis_nonidentifiability_1975}, or \emph{vice versa}, is a source of debate, see e.g. \cite{slud_nonparametric_1992}, \cite{fermanian_nonparametric_2003} and references therein. Leaving this interpretational controversy aside, they differ in an important statistical aspect: the \emph{crude} quantities in the mentioned papers of Peng and Fine do not suffer from the identifiability problem raised by \cite{tsiatis_nonidentifiability_1975}, that is, however, to be dealt with in our \emph{net} setting, which we believe can provide a valuable complement to the \emph{crude} one (see also \cite{jiang_estimating_2003}).

Another line of literature is devoted to semicompeting risks data. Translated to a censoring framework, this corresponds to $C$ being always observed for all individuals, even if $C > T$. In this case the latent margins (\emph{net} quantities) are more easily studied, as was also discussed in \cite{jiang_estimating_2003}. The article by \cite{li_quantile_2015} received most attention in the quantile regression literature and is relatively close to our model in that they also use a copula with unspecified association parameter, and linear marginal quantile models (but, as in our case, the linearity assumption for the censoring can be replaced). Yet, because of the additional information in their \emph{semi}competing risks setting, they can estimate it even while leaving the marginal distributions unspecified.
\bigbreak
Taking the final step to proper (mutual) censoring of $T$ and $C$ that is moreover not (conditionally) independent, the existing literature seems to reduce to a few articles only. The work by \cite{ji_analysis_2014} can be considered as the quantile counterpart of \cite{zheng_estimates_1995}, \cite{huang_regression_2008}, \cite{chen_semiparametric_2010} etc., working under the similar assumption of a fully known copula function coupling the marginal distributions. On the other hand, \cite{fan_partial_2018} leave the association parameter in their Archimedean copula models unspecified, but focus on identified \emph{sets} only, rather than point identification. Prior knowledge on the association strength can be used to narrow these sets for the regression coefficients, but no such sets are obtained for the copula parameter itself.

To the best of our knowledge, no quantile regression for the latent survival time $T|X$ subject to (conditionally) dependent censoring has been proposed where the association parameter is identifiable based on the data, which is the content of the current article.

\section{The model and notation} \label{sec:model}
The interest is in a model for (quantiles of) the survival time, that is subject to possibly dependent censoring. Throughout, we use the pair $(T,C)$ to denote the logarithms of the survival and censoring time, respectively, thanks to the transformation equivariance mentioned in Section \ref{sec:intro}. Only the first occurring one $Y = \min(T,C)$ is observed together with its status $\Delta = I(T \leq C)$, an indicator function taking on the value 1 if the survival time is observed and 0 if instead the censoring time is observed. We moreover assume the presence of a $(1+p)$-dimensional covariate $X = (1, \tilde{X})$ including an intercept entry, such that the observable data consist of triplets of $(Y,X, \Delta)$. For notational convenience we use one common covariate vector $X$ for $T$ and $C$, but one could also work with separate covariates $X$ for $T$ and $W$ for $C$, where the components of $X$ and $W$ can be the same, partially overlapping, or distinct. The range of $X$ is denoted $R$; it can always be decomposed as $R = \{1\} \times R_1 \times \dots \times R_p$, where each range $R_i$ of $\tilde{X}_i$ typically represents either the entire real line $\R$ or only a compact subset of it. Its values are always decomposed as $x = (1, \tilde{x})$. Throughout the paper, $F_Z(\cdot)$ denotes the distribution function of any random variable $Z$ with corresponding probability density $f_Z(\cdot)$ and quantile function $Q_Z(\cdot)$, where $Z$ can be just a single variable (e.g. $T$), a conditional version (e.g. $T|X$), a pair or triplet (e.g. $(Y,X, \Delta)$) or even a combination of both (e.g. $(T|X,C|X)$). As the conditioning is always done simultaneously on both variables in a pair, we slightly abuse notation and write e.g. $(T,C|X)$ rather than $(T|X, C|X)$ or $(T,C)|X$. Moreover, parentheses are typically omitted in subscripts, e.g. $F_{T,C|X}(\cdot, \cdot| \cdot)$. When conditioning on two variables, we simply write e.g. $F_{Y|X, \Delta}(\cdot| \cdot, \cdot)$. Finally, the subscripts are also used to indicate the parameters of the corresponding distributions using a semicolon, e.g. $f_{T|X; \theta_T}(\cdot| \cdot)$ refers to the density of $T|X$ with parameter vector $\theta_T$.

\subsection{Model formulation}\label{sec:model:form}
As explained in Section \ref{sec:intro}, we propose a fully parametric model for $(T|X, C|X) = (T,C|X)$, using the copula approach as in the work of \cite{czado_dependent_2023}, i.e. we assume 
\begin{equation}\label{eq:model:joint}
	F_{T,C|X}(t,c|x) = \cpl_\theta(F_{T|X; \theta_T}(t|x), F_{C|X; \theta_C}(c|x))
\end{equation}
for a parametric copula function $\cpl_\theta(\cdot, \cdot)$ determined by a parameter $\theta$ that is assumed to be independent of the covariate $X$, and ranging over a parameter space $\Theta$. Similarly, the parametric marginal distributions have parameter spaces $\Theta_T$ and $\Theta_C$, elements of which (typically vectors) are referred to by their lower-case variants. The specific choice of the distribution for $C$ can vary, as long as its parameters can be identified in the sense of assumption \ref{ass:idf:LC:asyidf} below. This includes, for example, the case where the logarithmic censoring time is normally distributed with homoscedastic error term $\epsilon_C$, independent of $X$:
\begin{equation}\label{eq:model:C}
	C = X^\top \alpha + \sigma_C \erc, \quad \mbox{with} \quad \erc \perp\!\!\!\perp X,
\end{equation}
with $\erc \sim \mathcal{N}(0,1)$. This will be our default choice in simulations and examples, but the theory is by no means restricted to this; even the linear specification of the mean is not necessary. For $T$, on the other hand, we impose a more specific structure, closely related to linear quantile regression. More specifically, let
\begin{equation}\label{eq:model:T}
	T = X^\top \beta + \sigma(X; \gamma) \ert, \quad \mbox{with} \quad \ert \perp\!\!\!\perp X,
\end{equation}
where the error term $\ert$ is independent of $X$, and depends on the parameters $\tilde{\phi}$ and $\phi$ explained below, and on some $\lambda \in (0,1)$.
The variance term $\sigma(X; \gamma)$ can be used to introduce heteroscedasticity in the model. It is assumed to be of a known form, but with possibly unknown parameter $\gamma$ (that is, however, identifiable within $\sigma(X; \gamma)$, cf. assumption \ref{ass:idf:gamidf} below); our default choice is $\sigma(X; \gamma) = \exp(X^\top \gamma)$. Clearly, whenever $\ert$ is such that its quantile of level $\lambda$ equals zero, the corresponding quantile of level $\lambda$ for the survival time,
\begin{equation}\label{eq:model:qu}
	Q_{T|X}(\lambda|x) = x^\top \beta,
\end{equation}
is compatible with linear quantile regression. In order to benefit from this connection, the distribution imposed on $\ert$ is constructed starting from an AL distribution with location parameter $0$ and scale parameter $1$, motivated by the equivalence between quantiles and maximum likelihood estimation under this AL distribution \cite{kotz_laplace_2001}. More specifically, denote the so-called check-loss function with parameter $\lambda$ by $\rho_{\lambda}(z) = z(\lambda - I(z < 0))$, using the indicator function $I(\cdot)$. Then the AL density with parameter $\lambda$ can be written as
\begin{equation*}
	f_{AL}(y | \lambda) = \lambda(1-\lambda) \exp(-\rho_{\lambda}(y))
\end{equation*}
and it is clear that log-likelihood maximisation in this framework (while fixing $\lambda$) comes down to minimising the check-loss function $\rho_{\lambda}(\cdot)$. This, in turn, is equivalent to finding the $\lambda$-th quantile \cite{koenker_regression_1978}, but this equivalence breaks down under the presence of censoring. \cite{bottai_laplace_2010} tried to correct for censoring in their AL-based quantile regression, but a letter to the editor by \cite{koenker_letter_2011} revealed the inconsistency of their estimator. The recent work of \cite{kreiss_VK_submitted} shows that this problem can be tackled by replacing the AL distribution by an extension. Any distribution in their extended family of so-called \emph{Enriched (asymmetric) Laplace} distributions (which we abbreviate EAL throughout) has a density function with parameters $(\tilde{\phi}, \phi, \lambda)$ that is given by
\begin{equation}\label{eq:dens:eal}
	f_{\eal}(y | \lambda) = \lambda (1-\lambda) \begin{cases}
		e^{-\lambda y} \|\phi\|^{-2} \left(\sum_{k = 0}^{m} \phi_k L_k (\lambda y)\right)^2  & y > 0 \\
		e^{-(\lambda-1) y} \|\tilde{\phi}\|^{-2} \left(\sum_{k = 0}^{\tilde{m}} \tilde{\phi}_k L_k ((\lambda-1) y)\right)^2 & y \leq 0,
	\end{cases}
\end{equation}
where (throughout the paper) $\| \cdot \|$ denotes Euclidean norm, and $\phi = (1, \phi_1, \dots, \phi_m)$ and similarly $\tilde\phi = (1, \tilde\phi_1, \dots, \tilde\phi_{\tilde{m}})$, with possibly $\tilde{m} \neq m$. The functions $L_k(\cdot)$ denote the Laguerre orthonormal polynomials of degree $k$,
\begin{equation}\label{eq:lag:coeff}
	L_k(x) = \sum_{j = 0}^{k} \binom{k}{j} \frac{(-1)^j}{j!} x^{j},
\end{equation}
that are orthogonal with respect to the Laguerre weight $e^{-x}$ and that are moreover normalised, such that they are the unique polynomial family for which $\int_{0}^{\infty} L_k(x) L_j(x) e^{-x} dx$ equals $1$ if $j=k$ and 0 otherwise. This orthogonality with respect to the appropriate weight function is a key feature in showing that the property of the $\lambda$-th quantile of the AL($\lambda$) distribution being zero, is preserved by any $\eal(\tilde{\phi}, \phi, \lambda)$ distribution (cf. Lemma 2.2 in \cite{kreiss_VK_submitted}), which implies that any distribution in this family is still admissible for $\ert$. Moreover, by increasing the degrees $m$ and $\tilde{m}$ of the maximal degree polynomial on the positive and negative real axis, respectively, one can approximate any continuous density in the Hellinger sense (cf. Section 4.1 in \cite{kreiss_VK_submitted}), so it is also sufficiently flexible to capture the full behaviour of $T$. (Note that the density in \eqref{eq:dens:eal} is \emph{a priori} not necessarily continuous at the origin, but for the estimation we impose this using a constraint on the parameters $\tilde{\phi}$ and $\phi$, cf. Section \ref{sec:parest} below.) The parameters $\tilde{\phi}$ and $\phi$ are referred to as \emph{weights} or \emph{coefficients} of the Laguerre polynomials; these terms are never to be interpreted in the strict sense of the actual Laguerre coefficients in \eqref{eq:lag:coeff} -- those are referred to as \emph{polynomial coefficients} -- or the Laguerre weight $e^{-x}$ with respect to which the latter are defined by their orthogonality. In conclusion, we assume distribution functions for $T|X$ of the form
\begin{equation}\label{eq:Tdistr}
	F_{T|X; \theta_{T}}(y | x) = F_{\eal(\tilde{\phi}, \phi, \lambda)}\left(\frac{y - x^\top \beta}{\sigma(x; \gamma)}\right),
\end{equation}
where the parameter $\theta_T$ consists of the tuple $(\beta, \lambda, \tilde{\phi}, \phi, \gamma)$ or $(\beta, \tilde{\phi}, \phi, \gamma)$ whose components may, in turn, contain vectors rather than scalars. Only $\lambda \in (0,1)$ is always a scalar, at least if present (cf. Remark \ref{rem:lampara}). Typically, we use a decomposition $(\beta_0, \tilde{\beta})$ for the intercept part and the one corresponding to $\tilde{X}$, respectively, and similarly for $\gamma$.

Finally, for further use below, we introduce two important functions associated to the copula distribution, namely its partial derivatives 
\begin{equation*}
	h_{C|T; \theta}(v | u) = \frac{\partial}{\partial u} \cpl_\theta(u,v)
	\qquad \mbox{and} \qquad
	h_{T|C; \theta}(u | v) = \frac{\partial}{\partial v} \cpl_\theta(u,v),
\end{equation*}
referred to as $h$-functions below. As shown in \cite{czado_dependent_2023}, coupling the margins with these $h$-functions yields the conditional marginal distributions, i.e.,
\begin{equation*}
	h_{T|C}(F_T(t)|F_C(c)) = F_{T|C}(t|c) \qquad \mbox{and} \qquad h_{C|T}(F_C(t)|F_T(c)) = F_{C|T}(c|t).
\end{equation*}
In particular, these can be used to express the joint density of $(Y,\Delta|X)$ as
\begin{multline}\label{eq:lik:contr}
	\left[f_{T|X; \theta_T}(y|x) \cdot \left(1 - h_{C|T; \theta}(F_{C|X; \theta_C}(y|x) | F_{T|X; \theta_T}(y|x))\right)\right]^{\delta} \cdot \\
	\left[f_{C|X; \theta_C}(y|x) \cdot \left(1 - h_{T|C; \theta}(F_{T|X; \theta_T}(y|x) | F_{C|X; \theta_C}(y|x))\right)\right]^{1 - \delta},
\end{multline}
where the expression on the first line ($\delta = 1$) corresponds to an uncensored observation $(y,x,1)$ of $(Y,X,\Delta)$ and the one on the second line ($\delta = 0$) to a censored one (see also formula (7) in \cite{czado_dependent_2023}).

\begin{Remark}\label{rem:lampara}
	At this point, it may still be unclear why $\lambda$ is in the set of parameter vectors (and only in some cases) rather than fixed. Typically, the quantile level of interest -- for which often linearity is imposed -- is taken as a fixed value, and a model is associated to each such choice of $\lambda$. In Section \ref{sec:model:multiqu} below we explain why including it as a model parameter sometimes actually makes sense. One of the consequences is that only one, global model is to be imposed rather than a separate one for every quantile level of interest.
\end{Remark}

\begin{Remark}\label{rem:model}
	\begin{enumerate}[(i)]
		\item \label{rem:model:distr.reg} Whether the proposed model still fits within the quantile regression framework, may be debatable: in determining the quantiles of $T$, we first estimate all parameters in model (\ref{eq:model:joint}) - (\ref{eq:model:T}). Therefore, rather than just estimating some of the survival time $T$'s quantiles, regression is actually done on its entire distribution. Still, it is clear that the model -- that we actually started developing with a purely quantile-oriented perspective in mind -- is intimately connected to quantile regression, which distinguishes it from \enquote*{real} distributional regression in survival analysis as in \cite{delgado_distribution_2022}, for instance.
		\item \label{rem:model:tautrunc} This (parametrically) distributional touch comes with the advantage of being able to determine the quantile of \emph{any} level in $(0,1)$, whereas in previous literature truncation to a (significantly) smaller subset was commonly adopted to overcome identifiability issues \cite{peng_survival_2008,ji_analysis_2014,li_quantile_2015}. Some more advantages will become clear later on.
	\end{enumerate}
\end{Remark}

\subsection{Parameter status of the linear quantile level $\lambda$}\label{sec:model:multiqu}
Since our model is motivated by an underlying linear quantile model for $T$ on some level $\lambda$, one question that arises is how to generalise this model to multiple quantile levels. Is it possible, and necessary, to impose model (\ref{eq:model:joint}) - (\ref{eq:model:T}) on any quantile level $\lambda$ of interest? Or does information on one level suffice to derive the other quantiles, and can this level even be selected based on the data, without specifying it beforehand?

\subsubsection{Identifiability of the regression coefficient $\beta$: intuition}\label{sec:model:multiqu:intuition}
To answer these questions, we first gain some insight in the relation between the variance term and the shape of the quantile curves. We assume that for some $\epsilon_\lambda$ independent of $X$, with strictly increasing distribution function and $\lambda$-th quantile equal to zero (not necessarily EAL distributed),
\begin{equation}\label{eq:multiqu:Tmodel}
	T = X^\top \beta + \sigma(X;\gamma) \epsilon_\lambda,
\end{equation}
such that \eqref{eq:model:qu} holds. Using some basic computations, any other quantile can then be written as
\begin{equation}\label{eq:qushift}
	Q_{T|X}(p|x) = Q_{T|X}(\lambda|x) + \sigma(x; \gamma) Q_{\epsilon_\lambda}(p) = x^\top \beta + \sigma(x; \gamma) Q_{\epsilon_\lambda}(p),
\end{equation}
which reduces to \eqref{eq:model:qu} only for $p = \lambda$. We make a case distinction based on the formula for $\sigma(X; \gamma)$. For the purpose of illustration only, we restrict ourselves to the case where $\tilde{X}$ is one-dimensional and $\beta = (\beta_0,\beta_1)$ and $\gamma = (\gamma_0, \gamma_1)$. In Supplement S1.1, the visual reasoning of the present section is mathematically substantiated in the more general case with $\tilde{X} \in \R^p$.	
\bigbreak
\noindent
\textbf{Case 1: $\sigma(X; \gamma) = \gamma_0 + \gamma_1\tilde{X}$ is linear.} Rearranging (\ref{eq:qushift}), the $p$-th quantile can be written as 
\begin{equation*}
	Q_{T|X}(p|x) = \left[ \beta_0 + \gamma_0 Q_{\epsilon_\lambda}(p) \right] + \tilde{x}\left[\beta_1 + \gamma_1 Q_{\epsilon_\lambda}(p)\right]
\end{equation*}
for any $p$, which shows that all quantile curves in this case are linear. The unicity of the intercept and slope is determined by whether or not $\gamma_0$ or $\gamma_1$ is zero; all cases are illustrated in Figure \ref{fig:qucurves}. For each of those cases, $\sigma > 0$ is assumed (cf. assumption \ref{ass:idf:gamidf} below), such that indeed $Q_{T|X}(p|x)$ is increasing in $p$ as illustrated. This can be guaranteed by imposing some extra conditions on $\gamma$, depending on the range of $\tilde{X}$. In each case, the intercepts (slopes) can be identified without specifying $\lambda$ if and only if they are the same for all curves in the corresponding panel. In particular, whenever at least one component of $\gamma$ is nonzero, there is always more than one linear quantile curve satisfying the model assumptions.

\begin{figure}[!h]
	\centering
	\subfigure{
		\begin{tikzpicture}[scale = 0.5, every node/.style={scale=0.7}]
			\filldraw[fill=black, draw=black] (3,7.5) circle (0cm) node[anchor=south] {\Large{$\boldsymbol{\gamma_1 = 0}$}};
			\draw[very thick,->] (0,0) -- (6,0) node[anchor=north, yshift = -5pt] {\Large{$\tilde{X}$}};
			\draw[very thick,->] (0,0) -- (0,6) node[anchor=south] {\Large{$Q_{T|X}(p|x)$}};
			
			\draw[very thick,-] (-0.9,0.6) -- (5.5,3) node[anchor=north, yshift = -5pt] {\Large{$p_1$}};
			\draw[very thick,-] (-0.9,1.3) -- (5.5,3.7) node[anchor=south] {\Large{$p_2 = \lambda$}};
			\draw[very thick,-] (-0.9,2.9) -- (5.5,5.3) node[anchor=south east] {\Large{$p_3$}};
			
			\filldraw[fill=black, draw=black] (0,0.95) circle (0.125cm) node[anchor=north west] {\Large{$\beta_0(p_1)$}};
			\filldraw[fill=black, draw=black] (0,1.65) circle (0.125cm) node[anchor=south east] {\Large{$\beta_0(p_2)$}};
			\filldraw[fill=black, draw=black] (0,3.25) circle (0.125cm) node[anchor=south east] {\Large{$\beta_0(p_3)$}};

			\filldraw[fill=black, draw=black] (11,7.5) circle (0cm) node[anchor=south] {\Large{$\boldsymbol{\gamma_0 = 0}$}};
			\draw[very thick,->] (8,0) -- (14,0) node[anchor=north, yshift = -5pt] {\Large{$\tilde{X}$}};
			\draw[very thick,->] (8,0) -- (8,6) node[anchor=south] {\Large{$Q_{T|X}(p|x)$}};
			
			\draw[very thick,-] (7.1,1.6) -- (13.5,2.2) node[anchor=north, yshift = -5pt] {\Large{$p_1$}};
			\draw[very thick,-] (7.1,1.3) -- (13.5,3.7) node[anchor=south east, xshift = 3pt, yshift = -3pt] {\Large{$p_2 = \lambda$}};
			\draw[very thick,-] (7.1,0.95) -- (13.5,5.9) node[anchor=south east] {\Large{$p_3$}};
			
			\filldraw[fill=black, draw=black] (8,1.65) circle (0.125cm) node[anchor=south east] {\Large{$\beta_0$}};
			
			\filldraw[fill=black, draw=black] (18.5,7.5) circle (0cm) node[anchor=south] {\Large{$\boldsymbol{\gamma_0 \neq 0 \neq \gamma_1}$}};
			\draw[very thick,->] (16,0) -- (22,0) node[anchor=north, yshift = -5pt] {\Large{$\tilde{X}$}};
			\draw[very thick,->] (16,0) -- (16,6) node[anchor=south] {\Large{$Q_{T|X}(p|x)$}};
			
			\draw[very thick,-] (15.1,1) -- (21.5,1.6) node[anchor=north, yshift = -5pt] {\Large{$p_1$}};
			\draw[very thick,-] (15.1,1.3) -- (21.5,4) node[anchor=south east] {\Large{$p_2 = \lambda$}};
			\draw[very thick,-] (15.1,2.5) -- (21.5,6.7) node[anchor=south east] {\Large{$p_3$}};
			
			\filldraw[fill=black, draw=black] (16,1.1) circle (0.125cm) node[anchor=north west] {\Large{$\beta_0(p_1)$}};
			\filldraw[fill=black, draw=black] (16,1.65) circle (0.125cm) node[anchor=south east] {\Large{$\beta_0(p_2)$}};
			\filldraw[fill=black, draw=black] (16,3.1) circle (0.125cm) node[anchor=south east] {\Large{$\beta_0(p_3)$}};
			
			\filldraw[fill=black, draw=black] (27,7.5) circle (0cm) node[anchor=south] {\Large{\textbf{nonlinear} $\boldsymbol{\sigma(\cdot)}$}};
			\draw[very thick,->] (24,0) -- (30,0) node[anchor=north, yshift = -5pt] {\Large{$\tilde{X}$}};
			\draw[very thick,->] (24,0) -- (24,6) node[anchor=south] {\Large{$Q_{T|X}(p|x)$}};
			
			\draw[very thick,dashed] (23.1,0.9) .. controls (28,-1) and (27,3) .. (29.5,2.7) node[anchor=north, yshift = -5pt] {\Large{$p_1$}};
			\draw[very thick,-] (23.1,1.3) -- (29.5,3.3) node[anchor=south east, xshift = 10pt] {\Large{$p_2 = \lambda$}};
			\draw[very thick,dashed] (23.1,2.5) .. controls (26,5.5) and (28,2.5) .. (29.5,5.8) node[anchor=south east] {\Large{$p_3$}};
			
			\filldraw[fill=black, draw=black] (24,0.6) circle (0.125cm) node[anchor=south west] {\Large{$\beta_0(p_1)$}};
			\filldraw[fill=black, draw=black] (24,1.65) circle (0.125cm) node[anchor=south east] {\Large{$\beta_0(p_2)$}};
			\filldraw[fill=black, draw=black] (24,3.25) circle (0.125cm) node[anchor=south east] {\Large{$\beta_0(p_3)$}};
		\end{tikzpicture}
	}
	\caption{Shape of the quantile curves depending on the form of $\sigma(X; \gamma)$, for Case 1 with $\gamma_1 = 0$, i.e. homoscedasticity (first panel), for $\gamma_0 = 0$ (second), $\gamma_0 \neq 0$ and $\gamma_1 \neq 0$ (third); and the nonlinear-$\sigma(\cdot)$ Case 2 (fourth panel).}
	\label{fig:qucurves}
\end{figure}
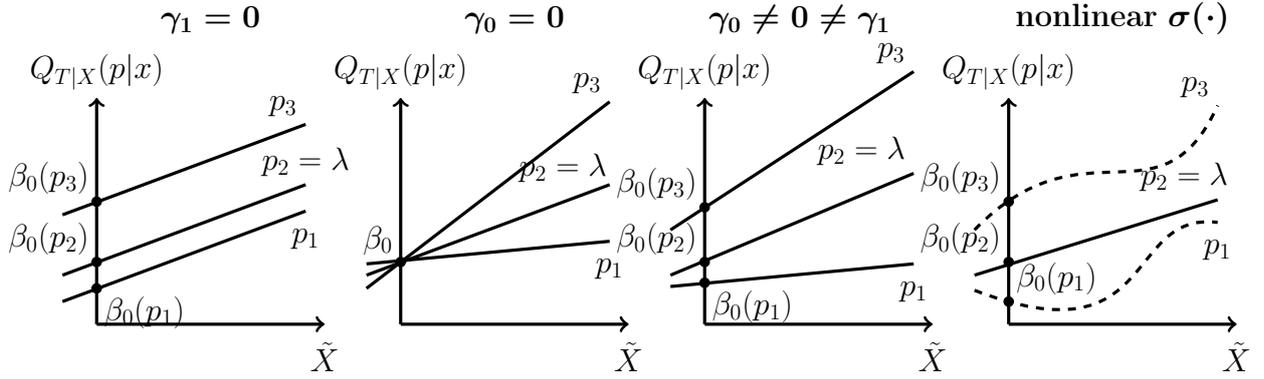

\noindent
\textbf{Case 2: $\sigma(X;\gamma)$ is nonlinear.} Equation (\ref{eq:qushift}) now cannot be rearranged in terms of a linear form, so only one $\lambda$ corresponds to a linear quantile curve (Figure \ref{fig:qucurves}). Therefore, even without fixing $\lambda$, a sufficient amount of data should enable identification of $\lambda$ as well as the regression curve.

\subsubsection{Implications: the parameter status of $\lambda$}\label{sec:model:multiqu:parastatus:lam}
The above discussion motivates a case distinction for $\lambda$, that is considered either as a model variable, or as a fixed parameter according to the variance type. It moreover has some interesting consequences for modelling several quantile levels of the same distribution; those are discussed in Section \ref{sec:model:multiqu:num}.
\bigbreak
Suppose that the variance is a general, nonlinear function of $X$. Then $\lambda$ need not -- and \emph{should} not -- be specified beforehand, but has to be a model parameter, since it is identifiable from the data. Indeed, if model (\ref{eq:multiqu:Tmodel}) is valid for some $\lambda$, imposing the same model (\ref{eq:multiqu:Tmodel}) for any other level $p \neq \lambda$ (with a possibly different set of parameters) would be in contradiction with the unicity of the linear quantile curve (Section \ref{sec:model:multiqu:intuition}, Case 2). Even when considering only one quantile level, knowing at precisely which unique level linearity occurs, seems very unlikely.

Consider on the other hand the case where the variance function is linear, such that \emph{all} quantile curves are linear (Section \ref{sec:model:multiqu:intuition}). Choosing a different $\lambda$ can always be made up for by adapting $\beta$; model (\ref{eq:multiqu:Tmodel}) can be identifiable only if $\lambda$ is fixed. On the other hand, assuming that the model is true for one value of $\lambda$ is equivalent to assuming that it holds on \emph{all} levels (not necessarily with the same parameter values). Thus, fixing $\lambda$ is necessary, but it also cannot be done in a wrong way. (Although model performance in practice is apparently not as independent of the chosen $\lambda$ as it is theoretically. This is discussed later on in the context of some simulation results in Section \ref{sec:sim}).
\bigbreak
In conclusion, if $\sigma(\cdot)$ is nonlinear, $\lambda$ should be introduced into the model as a variable parameter, whereas it should be included as a known constant in the linear case. These situations are in the remainder of the paper referred to as the \enquote*{$\lambda$ variable} and \enquote*{$\lambda$ fixed} case.
\begin{Remark}
	The present discussion assumed a general $\epsilon_\lambda$ that is constrained only to have $\lambda$-th quantile equal to zero, whereas we always fit models where the error term $\ert$ of \eqref{eq:model:T} follows an EAL distribution. Under fixed Laguerre degrees, imposing this EAL family removes some of the flexibility and the corresponding nonidentifiability. However, when fitting the model, these degrees are not specified beforehand, and the mentioned EAL flexibility for increasing degrees could raise the same identifiability issues. We therefore adhere to the parameter statuses discussed here (even for the identifiability proof, where the degrees $\tilde{m}$ and $m$ are still fixed).
\end{Remark}

\subsubsection{Implications: (numerical) multiple quantile modelling}\label{sec:model:multiqu:num} 
In the case of multiple quantile levels of interest, imposing model (\ref{eq:model:joint})-(\ref{eq:model:T}) for all those distinct levels comes at a large computational cost: the model parameters need to be estimated on each level $\lambda$ separately. This does not only entail time-consuming optimisation procedures, but also the undesirable possibility of crossing quantile curves. Fortunately, equation (\ref{eq:qushift}) enables reducing the numerical procedure to one optimisation only. The variable $\lambda$ optimisation yields an optimal value $\lambda_0$ for the parameter $\lambda$ as well as a corresponding set of EAL parameters. These can be used to estimate the linear $\lambda_0$-quantile, and on any other level $p$ one only has to compute the shift that depends on $p$, $\lambda_0$ and the optimal EAL parameters. On the other hand, when $\lambda$ needs to be fixed, one can just pick any quantile level, e.g. the median (cf. Section \ref{sec:model:multiqu:parastatus:lam}), to perform numerical optimisation and then once again apply the shift approach of (\ref{eq:qushift}) to obtain the other quantiles.

\section{Identifiability}\label{sec:idf}
Parameter estimation is done based on maximum likelihood estimation. We should thus first verify that the model is identifiable in the sense of necessary equality of two sets of parameter vectors whenever they yield the same contributions to the likelihood. More precisely, suppose that the observed data $\datan = \{(y_i, x_i,\delta_i), i = 1, \dots, n\}$ consist of $n$ observed tuples of $(Y,X,\Delta)$. Write $\pi = (\theta, \theta_T, \theta_C) \in \Pi = \Theta \times \Theta_T \times \Theta_C$ for the entire parameter vector, containing all components for the copula and both margins. We assume that $\lambda$ is omitted from the parameter vector $\theta_T = (\beta, \lambda, \tilde{\phi}, \phi, \gamma)$ in those cases for which it cannot be identifiable (cf. \emph{supra}). Denote the joint log-likelihood for the parameter $\pi$ with respect to $\datan$ by $\ell(\datan; \pi)$; in view of \eqref{eq:lik:contr}, it is of the form
\begin{eqnarray}\label{eq:loglik}
	\ell(\datan; \pi) & = & \sum_{i = 1}^n \log \left[f_{Y,\Delta | X; \pi}(y_i, \delta_i|x_i)\right]\\
	& = & \begin{multlined}[t]
		\sum_{\delta_i = 1} \log \Big[f_{T|X; \theta_T}(y_i|x_i) \cdot \left(1 - h_{C|T; \theta}(F_{C|X; \theta_C}(y_i|x_i) | F_{T|X; \theta_T}(y_i|x_i))\right)\Big] + \\
		\sum_{\delta_i = 0} \log \Big[f_{C|X; \theta_C}(y_i|x_i) \cdot \left(1 - h_{T|C; \theta}(F_{T|X; \theta_T}(y_i|x_i) | F_{C|X; \theta_C}(y_i|x_i))\right)\Big].
	\end{multlined} \nonumber
\end{eqnarray}
For the log-likelihood contribution of just one tuple $(y, x,\delta)$ rather than the entire collection of data $\datan$, we use the notation $\ell(y,x,\delta; \pi)$, which is nothing but $\ell(\mathcal{D}_1; \pi)$ in the case where $\mathcal{D}_1 = \{(y,x, \delta)\}$; it equals the logarithm of \eqref{eq:lik:contr}.

To show now that the joint model for $(T,C|X)$ with observations $(Y,X,\Delta)$ is identifiable as formally stated in Theorem \ref{thm:idf} and \ref{thm:idf:clayton} below, a reasoning similar to the one in \cite{czado_dependent_2023} can be applied. Not all sufficient conditions for their main result (Theorem 1) are satisfied: knowing that the limiting behaviour of $T$ is the same for two sets of parameters does not suffice to deduce equality of all parameters for $T$, but only part of them. However, the ideas in their proof can be modified using an \emph{ad hoc} method. Whenever all parameters in the censoring distribution are identifiable, this information can be used to first identify the copula parameter and only afterwards identify the remaining parameters of $T$'s distribution. A drawback of this \emph{ad hoc} proof is that it uses the explicit expression of the copula function and hence has to be repeated for each different copula function separately. Moreover, this proof method rules out using the Gaussian copula, lacking an explicit distribution function.

All assumptions for Theorem \ref{thm:idf} and \ref{thm:idf:clayton} are listed in Section \ref{sec:idf:ass}. Since some of them are rather technical, a discussion is included afterwards (Section \ref{sec:idf:intuition}). Some auxiliary lemmata to verify the assumptions are provided in Appendix \ref{app:idf:auxlem}, while Appendix \ref{app:proof:idfthm} contains the proof of the identifiability theorems.

\subsection{Assumptions and identifiability statements}\label{sec:idf:ass}
In order to formulate the identifiability assumptions, we first introduce some extra notation. Given any covariate value $X = x$, denote the set of outcomes $y$ with positive observation or censoring probability conditional on $X = x$ and $Y = y$ by $U_T(x)$ and $U_C(x)$, respectively, i.e.
\begin{equation*}
	U_T(x) = \{y | P(\Delta = 1 | Y = y, X = x) > 0\}, \quad \text{ and } \quad U_C(x) = \{y | P(\Delta = 0 |Y = y, X = x) > 0\}.
\end{equation*}
Note that when $U_T(x) = \R = U_C(x)$ for all $x \in \R$, this implies the assumptions $P(\Delta = 1 | Y,X) > 0$ and $P(\Delta = 0 | Y,X) > 0$. We deviate from these more standard forms to also allow for the case where the support of $C$ has a finite upper bound (the one for $T$ always spans the entire real line by our EAL model assumption). In that case, $U_T(x) \subsetneq \R$ and similarly $U_C(x)$ is strictly contained in $\R$, for any value $X = x$. Next, denote the support of $C|X$ by $\mbox{supp} f_{C|X; \theta_C} = [M_l(X, \theta_C), M_r(X, \theta_C)]$. Both bounds may attain infinite values for some or all of the $x$; typically, $M_l(x, \theta_C) = - \infty$ for any $x$, meaning that, regardless of the covariate value, censoring may occur from time zero onwards (recall the logarithmic notation for $C$.) For convenience, we implicitly assume in \ref{ass:idf:LC} below that $\mbox{supp} f_{C|X; \theta_C} = [M_l, M_r]$ depends on neither $\theta_C$ nor the value $X = x$, but this can be relaxed (cf. Remark \ref{rem:dep:XnotPara}). Finally, define the set
\begin{equation*}
	S(\beta) = \{x | x^\top \beta < M_r\},
\end{equation*}
where $\beta$ is the part of $\theta_T$ corresponding to the regression coefficients for $T$. The general list of identifiability assumptions is the following:

\begin{enumerate}[label = (A\arabic{enumi})]
	\item \label{ass:idf:posdef} The covariate $X = (1,\tilde{X})$ is such that $\Var (\tilde{X})$ is positive definite.
	\item \label{ass:idf:gamidf} The variance term $\sigma(X; \gamma)$ is such that $\sigma(x; \gamma) > 0$ for all $x$, and such that $\sigma(x; \gamma_1) = \sigma(x; \gamma_2), \forall x \in R$ implies that $\gamma_1 = \gamma_2$.
	\item \label{ass:idf:LC} For $L_C$ equal to at least one of the endpoints $\{M_l, M_r\}$ of $C|X$'s support (possibly infinite constants), the following three conditions are simultaneously satisfied:
	\begin{enumerate}[label = (A\arabic{enumi}.\alph{enumii})]
		\item \label{ass:idf:LC:seq} For any $x \in R$, $U_C(x)$ contains an interval $I_x$ of the form $I_x = (M_l, A_x)$, in case $L_C = M_l$, or of the form $I_x = (A_x, M_r)$, when $L_C = M_r$ (for some $A_x \in \R$).
		\item \label{ass:idf:LC:hfun} For any $x \in R$, for any copula parameter $\theta \in \Theta$ and $(\theta_T, \theta_C) \in \Theta_T \times \Theta_C$:
		\begin{equation*}
			\lim\limits_{y \to L_C} h_{T|C; \theta}(F_{T|X; \theta_T}(y|x)| F_{C|X; \theta_C}(y|x)) = 0.
		\end{equation*}
		\item \label{ass:idf:LC:asyidf} The parameters of $C$ are identifiable from the limiting behaviour at $L_C$, i.e. for any $\theta_{C_1}, \theta_{C_2} \in \Theta_C$,
		\begin{equation*}
			\lim\limits_{y \to L_C} \frac{f_{C|X; \theta_{C_1}}(y|x)}{f_{C|X; \theta_{C_2}}(y|x)} = 1, \forall x \in R \Iff \theta_{C_1} = \theta_{C_2}.
		\end{equation*}
	\end{enumerate}
	
	\item \label{ass:idf:LT} The following two conditions are simultaneously satisfied:
	\begin{enumerate}[label = (A\arabic{enumi}.\alph{enumii})]
		\item \label{ass:idf:LT:seq} For any $x \in R$, $U_T(x)$ contains an interval $J_x = (-\infty, B_x)$ for some $B_x \in \R$.
		\item \label{ass:idf:LT:hfun} For any $x \in R$, for any copula parameter $\theta \in \Theta$ and $(\theta_T, \theta_C) \in \Theta_T \times \Theta_C$:
		\begin{equation*}
			\lim\limits_{y \to - \infty} h_{C|T; \theta}(F_{C|X; \theta_C}(y|x)| F_{T|X; \theta_T}(y|x)) = 0.
		\end{equation*}
	\end{enumerate}
	
	\item \label{ass:idf:coppara} For any two pairs of candidates $(\theta_i, \theta_{T_i})$ with $i = 1,2$, there holds: whenever the parameter $\theta_{C}$ is identified and the below identity holds for any $x$, any corresponding $y \in U_C(x)$, then $\theta_1 = \theta_2$ is implied:
	\begin{equation*}
		1 = \frac{h_{T|C; \theta_1}(F_{T|X; \theta_{T_1}}(y|x) | F_{C|X; \theta_{C}}(y|x))}{h_{T|C; \theta_2}(F_{T|X; \theta_{T_2}}(y|x) | F_{C|X; \theta_{C}}(y|x))}.
	\end{equation*}
	\item \label{ass:idf:S0} For any candidate value $\beta$, there holds: whenever $x \in S(\beta)$, the point $x^\top \beta$ belongs to $U_C(x)$.
	\item \label{ass:idf:S1} For the true value $\beta^0$, there moreover exist some $x^0 \in S(\beta^0)$ and $y_L < y_R$ such that $U_C(x^0)$ contains the entire interval $I = (y_L, y_R)$ and $(x^0)^\top \beta^0$ is an internal point of $I$.
	\item \label{ass:idf:S2} For any pair of candidate values $\beta_i$, $i = 1,2$, the intersection $S(\beta_1) \cap S(\beta_2)$ is nonempty and moreover $\Var( \tilde{X} | X \in S(\beta_1) \cap S(\beta_2))$ is positive definite.
\end{enumerate}
All limits in assumptions \ref{ass:idf:LC}-\ref{ass:idf:LT} are to be interpreted as one-sided limits $y \to L^-$ (for $L = M_r$) or $y \to L^+$ (for $L \in \{M_l, -\infty \}$) from the negative or positive side, respectively. Moreover, for the variable $\lambda$ case, \ref{ass:idf:S2} may be omitted; when $\lambda$ is fixed, \ref{ass:idf:S2} is to be included, but \ref{ass:idf:S0}-\ref{ass:idf:S1} can be dropped.

\begin{Theorem}\label{thm:idf}
	Suppose that assumptions \ref{ass:idf:posdef}-\ref{ass:idf:S2} are satisfied and suppose that two candidate sets of parameters $\pi_1 = (\theta_1, \theta_{T_1}, \theta_{C_1})$ and $\pi_2 = (\theta_2, \theta_{T_2}, \theta_{C_2})$ yield the same contribution $\ell(y,x,\delta; \pi_1) = \ell(y,x,\delta; \pi_2)$ to the likelihood \eqref{eq:loglik} for any value $(y, x, \delta)$ of $(Y,X,\Delta)$. Then $\pi_1 = \pi_2$.
\end{Theorem}
As discussed in Appendix \ref{app:idf:auxlem}, assumptions \ref{ass:idf:posdef}-\ref{ass:idf:S2} cannot all be satisfied simultaneously in the case of a Clayton copula. In that case, we have the list of alternative assumptions
\begin{enumerate}[label = (B\arabic{enumi})]
	\item \label{ass:idf:B1} The Clayton copula parameter $\theta$ is strictly positive.
	\item \label{ass:idf:B2} The margins are such that $\lim\limits_{y \to - \infty} \frac{\log F_{C|X; \theta_C}(y|x)}{\log F_{T|X; \theta_T}(y|x)} = \infty$.
	\item \label{ass:idf:B3} For any $x \in R$, $U_C(x)$ contains an interval of the form $I_x = (-\infty, A_x)$.
	\item \label{ass:idf:B4} The parameters of $C$ are identifiable from its behaviour at $- \infty$, i.e. for any $\theta_{C_1}, \theta_{C_2} \in \Theta_C$,
	\begin{equation*}
		\lim\limits_{y \to -\infty} \frac{f_{C|X; \theta_{C_1}}(y|x)}{f_{C|X; \theta_{C_2}}(y|x)} = 1, \forall x \in R \Iff \theta_{C_1} = \theta_{C_2}.
	\end{equation*}
\end{enumerate}
They can be used to formulate the alternative (Clayton) identifiability theorem below.
\begin{Theorem}\label{thm:idf:clayton}
	Suppose that $\cpl_\theta(\cdot, \cdot)$ is a Clayton copula and that assumptions \ref{ass:idf:posdef}-\ref{ass:idf:gamidf}, \ref{ass:idf:S0}-\ref{ass:idf:S2} and the alternative assumptions \ref{ass:idf:B1}-\ref{ass:idf:B4} are satisfied. Then whenever two candidate sets of parameters $\pi_1 = (\theta_1, \theta_{T_1}, \theta_{C_1})$ and $\pi_2 = (\theta_2, \theta_{T_2}, \theta_{C_2})$ yield the same contribution $\ell(y,x,\delta; \pi_1) = \ell(y,x,\delta; \pi_2)$ to the likelihood \eqref{eq:loglik} for any value $(y, x, \delta)$ of $(Y,X,\Delta)$, this implies $\pi_1 = \pi_2$.
\end{Theorem}
Appendix \ref{app:idf:auxlem} contains a number of auxiliary lemmata that can be used to assess the validity of some of the identifiability assumptions. It moreover clarifies how \ref{ass:idf:B1}-\ref{ass:idf:B4} arise as an alternative to \ref{ass:idf:LC}-\ref{ass:idf:coppara}.

\subsection{Discussion on the assumptions}\label{sec:idf:intuition}
All identifiability assumptions ensure that enough information is available, amongst others by guaranteeing that the regions with positive observation probability are sufficiently rich (cf. \ref{ass:idf:LT}), and so are those with positive censoring probability (cf. \ref{ass:idf:LC}, \ref{ass:idf:S0}-\ref{ass:idf:S2}).
Despite their rather technical formulation, their content is quite natural. For example, assumption \ref{ass:idf:S2} essentially states that the proportion of the regression line $Y = X^\top \beta$ lying beyond the observable region cannot be too large, and the available part cannot contain too little information. Moreover, some conditions have been included in a more concise or simple form, but may actually be generalised or weakened a bit. Some possible generalisations are discussed in the following remarks. Conversely, assumption \ref{ass:idf:LC} has been formulated rather generally, but, as is apparent from the auxiliary lemmata (Appendix \ref{app:idf:auxlem}), we can in most cases just take $L_C = -\infty$. Note that together with \ref{ass:idf:LT}, this translates to both censoring and death occurring with nonzero probability right from time zero onwards, which doesn't seem too unnatural of an assumption.

\begin{Remark}[Generalisation of \ref{ass:idf:S1}]
	As can be seen from the proof (Appendix \ref{app:proof:idfthm}), it is actually sufficient to find an interval $I_- \subset U_C(x_-)$ below $(x_-)^\top \beta^0$ for some $x_-$ and another interval $I_+ \subset U_C(x_+)$ above $(x_+)^\top \beta^0$ for a possibly different value $x_+$. There is hence no need for both intervals to correspond to the same point $x^0$, nor for either of them to contain their corresponding $(x^0)^\top \beta^0$-value.
\end{Remark}

\begin{Remark}[Generalisation of \ref{ass:idf:LC}: dependence on $X$]\label{rem:dep:XnotPara}
	\begin{enumerate}[(i)]
		\item \label{rem:dep:XnotPara:notPara} First, notice that the proof method in Appendix \ref{app:proof:idfthm:C} requires that $L_C$ be the same for both candidate values for $\theta_C$: the identification of $C$'s parameters involves taking a limit for $y \to L_C$ in the expression
		\begin{equation*}
			\frac{f_{C|X; \theta_{C_1}}(y|x)}{f_{C|X; \theta_{C_2}}(y|x)} \cdot \frac{1 - h_{T|C; \theta_1}(F_{T|X; \theta_{T_1}}(y|x) | F_{C|X; \theta_{C_1}}(y|x))}{1 - h_{T|C; \theta_2}(F_{T|X; \theta_{T_2}}(y|x) | F_{C|X; \theta_{C_2}}(y|x))}
		\end{equation*}
		of equation (\ref{eq:proof:gumbel:cens:llh}) and relies on both $h$-functions involved tending to zero. This is guaranteed in both the numerator and the denominator only when both have the same endpoint of the support; otherwise the limit reasoning using assumption \ref{ass:idf:LC:asyidf} breaks down.
		\item \label{rem:dep:XnotPara:X} Strictly speaking, dependence of $L_C$ on $x$ for this reasoning -- and elsewhere in the identifiability proof -- is allowed, if also in assumption \ref{ass:idf:LC:asyidf} $L_C$ can be replaced by $L_C(x)$, but a situation where the upper bound of $C|X$ depends on $X$ but not on $\theta_C$ does not seem very likely. (In that case, also the assumptions involving the sets $S(\beta)$ need to be treated with care, since the latter sets will depend on $M_r(X)$. Both the assumptions and the corresponding parts in the proof can be modified to make this work, but this version has not been included in the paper.)
		\item \label{rem:dep:XnotPara:suppC} Thanks to Lemma \ref{lem:hfun:lim}, we can usually work with $L_C = M_l = - \infty$, where independence on $\theta_C$ and even $X$ is clearly satisfied. Only in the Gumbel case, the upper bound of $C|X$'s support is used rather than its lower bound. Both \ref{ass:idf:LC} and Lemma \ref{lem:hfun:lim} are formulated in such a way that $\mbox{supp} f_{C|X; \theta_C} = [M_l, M_r]$ depends on neither $X$ nor $\theta_C$, and moreover has a finite upper bound. This can be further relaxed. In the Gumbel case, the lower bound $M_l(X, \theta_C)$ can actually be anything (in all other cases, this goes for the upper bound $M_r(X, \theta_C)$ instead). Moreover, for the upper bound, in view of part \ref{rem:dep:XnotPara:notPara} and \ref{rem:dep:XnotPara:X} of this remark, also $M_r(X)$ would be allowed, as long as for any two candidate parameter values for $C$, $M_r(X, \theta_{C_1}) = M_r(X, \theta_{C_2})$. Either of the following conditions guarantees this:
		\begin{enumerate}[label = (G\arabic{enumii})]
			\item \label{rem:GumbelSuppC:const} The upper bound $M_r(X, \theta_C) = M_r < \infty$ is a finite constant.
			\item \label{rem:GumbelSuppC:alternative} For any $x \in R$ and any parameter value $\theta_C$, the density of $C|X = x$ is such that $M_r(x, \theta_C) < \infty$ and
			\[\forall y < M_r(x, \theta_C): \quad f_{C|X; \theta_C}(y|x) > 0, \qquad \text{but} \qquad \lim\limits_{y \to M_r(x, \theta_C)^-} f_{C|X; \theta_C}(y|x) = 0, \]
			where $y \to M_r(x, \theta_C)^-$ denotes the limit for $y$ tending to $M_r(x, \theta_C)$ from the left.
		\end{enumerate}
		In Supplement S1.2.6 we show that, assuming only \ref{rem:GumbelSuppC:alternative}, the fact that two sets of candidate parameters lead to the same likelihood contributions automatically entails that the corresponding upper bounds $M_r(X,\theta_{C_1})$ and $M_r(X,\theta_{C_2})$ are the same. It seems reasonable to assume that in the case of an upper-bounded support, either $C$ is somehow \enquote*{truncated} by another (external) mechanism that does not depend on the precise value of $\theta_C$ (cf. assumption \ref{rem:GumbelSuppC:const}), or it is naturally bounded by this upper bound $M_r(X, \theta_C)$ -- now possibly dependent on $\theta_C$'s value -- and its density function continuously flattens out towards $M_r(X,\theta_C)$, gradually decreasing to zero rather than just dropping jumpwise to zero at $M_r(X,\theta_C)$ (cf. \ref{rem:GumbelSuppC:alternative}).
	\end{enumerate}
\end{Remark}
Still, the easiest case to satisfy the Gumbel assumptions is for $C|X$ to have a constantly upper-bounded support, independent of $\theta_C$, cf. \ref{rem:GumbelSuppC:const}. This in particular rules out using the \enquote*{default model} $C = x^\top \alpha + \epsilon_C$ for $C$ proposed in (\ref{eq:model:C}): even if $\epsilon_C$ has a constant, finite upper bound, the support of $C|X$ depends on $\alpha$. However, a natural example of a model that is very close to (\ref{eq:model:C}) but \emph{does} satisfy the assumption on the constant upper bound with respect to $\theta_C$ is discussed in Supplement S1.3. It furthermore serves the purpose of showing the validity of a number of assumptions in a natural situation, moreover for the seemingly most restrictive case of the Gumbel family.

\section{Asymptotics}\label{sec:asy}
The consistency and asymptotic normality of the estimator in a maximum likelihood framework are typically obtained using the theory of \cite{white_maximum_1982}. However, the standard regularity conditions given there are not all satisfied: the likelihood function is not differentiable (let alone continuously so) in all arguments because of the two-sided definition of the EAL density. On the other hand, differentiability only breaks down for points on the regression line $Y = X^\top \beta$, enabling application of results for slightly more general contexts as in e.g. \cite{newey_large_1994} or \cite{pakes_simulation_1989}; we apply the results of the latter. We start by introducing the notions involved and giving a set of assumptions under which consistency and asymptotic normality can be proved using Corollary 3.2 and Theorem 3.3, respectively, in \cite{pakes_simulation_1989}. Next, we explain how our estimator is defined, and in the subsequent sections, we discuss its consistency and asymptotic normality.

\subsection{Assumptions and notation}
Denote $q$ the total length of the parameter vector $\pi = (\theta, \theta_T, \theta_C)$, i.e. the dimension of the corresponding parameter space $\Pi = \Theta \times \Theta_T \times \Theta_C$. Next to the log-likelihood $\ell(\cdot; \pi)$ already introduced in \eqref{eq:loglik}, we also need its derivatives with respect to the components of $\pi$, i.e. the score functions $\psi_j(y,x,\delta; \pi) = \frac{\partial \ell}{\partial \pi_j} (y,x,\delta; \pi)$ for $j = 1, \dots, q$ (with obvious extension to the full set of data $\datan$). On the regression line $Y = X^\top \beta$, these are defined as the mean of the left and right derivative: the constraints imposed in Section \ref{sec:parest} below ensure continuity of the EAL density and hence of the complete log-likelihood $\ell(\cdot; \pi)$, but not differentiability. Dependence of both functions $\ell(\cdot)$ and $\psi(\cdot)$ on the parameter $\pi$ is sometimes suppressed, unless ambiguity could arise. For $j = 1, \dots, q$, we consider the classes of score functions
\begin{equation}\label{eq:def:funclass}
	\mathcal{F}_j = \{\psi_j(\cdot; \pi): \R \times R \times \{0,1\} \to \R : (y,x,\delta) \mapsto \frac{\partial \ell}{\partial \pi_j} (y, x, \delta) \ | \ \pi \in \Pi\},
\end{equation}
that can be further decomposed because of the different form of the likelihood for censored observations ($\delta=0$) versus uncensored ones ($\delta = 1$). Indeed, we can also write
\begin{equation}\label{eq:psi:division}
	\psi_j(y,x,\delta; \pi) = I(\delta = 0) \cdot \psi_j^0(y,x; \pi) + I(\delta = 1) \cdot \psi_j^1(y,x; \pi).
\end{equation}
For the parts $\psi_j^0(\cdot, \cdot)$ and $\psi_j^1(\cdot, \cdot)$ corresponding to the censored and observed cases, respectively, we introduce
\begin{equation}\label{eq:def:funclass:delta}
	\mathcal{F}_j^0 = \{\psi_j^0(\cdot; \pi): \R \times R \to \R : (y,x) \mapsto \frac{\partial \ell}{\partial \pi_j} (y,x, 0) \ | \ \pi \in \Pi\},
\end{equation}
and similarly for $\mathcal{F}_j^1$.
\bigbreak
In all assumptions below, $\pi_0$ denotes the true parameter vector. For expectations, $E_{Z}[\cdot]$ means that the expectation is taken with respect to the (true) distribution of $Z$. This subscript is explicitly included only if the variable differs from the default triple $(Y,X,\Delta)$, or to avoid ambiguity. Similarly, unless mentioned otherwise, $P$ denotes probability with respect to the distribution of $(Y,X, \Delta)$.

\subsubsection{Assumptions: consistency}
\begin{enumerate}[label = (C\arabic{enumi})]
	\item \label{ass:cons:Cpiecectu} The distribution of $C$ is such that $f_C(\cdot)$ is piecewise continuous.
	\item \label{ass:cons:min.loc.unique} For each $\delta > 0$, $ \inf \limits_{|| \pi - \pi_0 || > \delta} ||E[\psi(Y,X, \Delta; \pi)]|| > 0$.
	\item \label{ass:cons:compactpara} The parameter space $\Pi = \Theta \times \Theta_T \times \Theta_C$ is compact.
	\item \label{ass:cons:domfun} There exists a dominating function $d(y,x,\delta)$ such that $d(Y,X,\Delta)$ has a finite mean and $||\psi(y,x, \delta; \pi)|| \leq d(y,x,\delta)$ for all $\pi \in \Pi$. 
	\item \label{ass:cons:Leibniz} For any admissible parameter $\pi$, Leibniz integral rule can be applied on either side of the regression line, for function \eqref{eq:lik:contr} evaluated in both $\delta = 0$ and $\delta = 1$.
\end{enumerate}

\subsubsection{Assumptions: asymptotic normality}
On top of the consistency assumptions, assume that
\begin{enumerate}[label = (N\arabic{enumi})]
	\item \label{ass:norm:fullrank} $E[\psi(Y,X,\Delta; \pi)]$ has a derivative matrix $\Gamma$ at $\pi_0$ that has full rank.
	\item \label{ass:norm:intpoint} $\pi_0$ is an interior point of $\Pi = \Theta \times \Theta_T \times \Theta_C$.
	\item \label{ass:norm:Donsk} For each $j = 1, ..., q$, the function class $\mathcal{F}_j$ satisfies the $P$-Donsker property.
	\item \label{ass:norm:Fjfinsup} Any such class $\mathcal{F}_j$ moreover satisfies $\sup\limits_{f \in \mathcal{F}_j} |E[f(Y,X,\Delta)]| < \infty$.
	\item \label{ass:norm:L2Pctu} For $j = 1, \dots, q$, the function $\psi_j(\cdot)$ is $L_2(P)$-continuous at $\pi_0$, i.e. whenever $|| \pi - \pi_0 || \to 0$,
	\[E[ (\psi_j(Y,X,\Delta; \pi) - \psi_j(Y,X,\Delta; \pi_0))^2] \to 0.\]
\end{enumerate}

Sufficient conditions for assumption \ref{ass:norm:L2Pctu} are discussed in Remark S2.1 in Supplement S2.4. Whether assumption \ref{ass:norm:Donsk} is satisfied is not obvious in general. It is guaranteed, among others, by the set of assumptions (D1)-(D4) provided in Supplement S2.4.2 together with a proof. Finally, assumption \ref{ass:norm:deltamethod} below is necessary only to transfer normality of the parameter estimators to that of the quantiles (Corollary \ref{cor:deltamethod}).

\begin{enumerate}[label = (N\arabic{enumi}), resume]
	\item \label{ass:norm:deltamethod} The variance term $\sigma(X; \gamma)$ is differentiable with respect to all components of $\gamma$.
\end{enumerate}

\subsection{Definition of the estimators $\hat{\pi}_n$}
The usual definition of the estimator $\hat{\pi}_n$ in terms of zeros of the score function $\psi(\cdot) = (\psi_1(\cdot), \dots, \psi_q(\cdot))$ is problematic in our EAL context. The optimal parameter $\beta$ (part of $\pi$) now determines the location of the discontinuity of $\psi(\cdot)$ across the regression line, rather than its zeros. Yet, this approach can be extended to (roughly) finding a $\pi$ minimising the norm of the derivative rather than being a proper zero, choosing the function $G_n(\cdot)$ of Pakes and Pollard's Corollary 3.2 as the empirical (componentwise) mean of the gradient. To make this idea more precise, define the vector function $G_n(\mathcal{D}_n; \pi) = (G_{n,1}, \dots, G_{n,q})(\mathcal{D}_n; \pi)$ whose component functions for $j = 1, \dots, q$ are given by
\[G_{n,j}(\mathcal{D}_n; \pi) = \frac{1}{n} \sum_{i = 1}^{n} \psi_j(y_i, x_i, \delta_i; \pi)\]
given the set of observations $\mathcal{D}_n$ as above. Now, for a fixed set of data $\mathcal{D}_n$, in the following omitted from the notation, define the estimator $\hat{\pi}_n$ to be any parameter value in the parameter space $\Pi$ such that
\begin{equation}\label{eq:def:pihat}
	||G_n(\hat{\pi}_n) || \leq \inf\limits_{\pi \in \Pi} ||G_n(\pi)|| + \frac{1}{n}.
\end{equation}
Intuitively, one would want to define $\hat{\pi}_n = \arg \min ||G_n(\pi)||$, but its existence is not guaranteed in the current setting, as opposed to $\hat{\pi}_n$ as in \eqref{eq:def:pihat}. We elaborate a little on these assertions in Supplement S2.1.

\subsection{Consistency}
We prove consistency of the estimator $\hat{\pi}_n$ of the previous section by verifying the conditions of Corollary 3.2 in \cite{pakes_simulation_1989}. This is done under assumptions \ref{ass:cons:Cpiecectu}-\ref{ass:cons:Leibniz}, and $G$ is the theoretical mean of the score functions, counterpart of the empirical $G_n$. Supplement S2.2 contains the proof for the consistency theorem stated below.

\begin{Theorem}\label{thm:cons}
	Under assumptions \ref{ass:cons:Cpiecectu}-\ref{ass:cons:Leibniz}, the estimator $\hat{\pi}_n$ in (\ref{eq:def:pihat}) converges in probability to the true parameter $\pi_0$.
\end{Theorem}

\subsection{Asymptotic normality}
Given consistency of the estimator, we use Theorem 3.3 in \cite{pakes_simulation_1989} to prove the asymptotic normality stated in the next theorem; the proof is deferred to Supplement S2.4.
\begin{Theorem}\label{thm:asynorm}
	Under assumptions \ref{ass:cons:Cpiecectu}-\ref{ass:cons:Leibniz} and \ref{ass:norm:fullrank}-\ref{ass:norm:L2Pctu}, the estimator $\hat{\pi}_n$ in (\ref{eq:def:pihat}) satisfies
	\[\sqrt{n} (\hat{\pi}_n - \pi_0) \stackrel{d}{\to} \mathcal{N}_q\left(0, (\Gamma^\top \Gamma)^{-1} \Gamma^\top V \Gamma (\Gamma^\top \Gamma)^{-1} \right), \]
	where $\stackrel{d}{\to}$ denotes convergence in distribution, $\Gamma$ is the derivative matrix mentioned in \ref{ass:norm:fullrank}, and $V$ is the asymptotic variance matrix such that $\sqrt{n} G_n(\pi_0) \stackrel{d}{\to} \mathcal{N}_q(0,V)$, arising from the central limit theorem, i.e. for any $i,j = 1, \dots, q$,
	\[V_{ij} = E\left[\frac{\partial \ell}{\partial \pi_i}\Big(Y, X, \Delta; \pi_0 \Big) \frac{\partial \ell}{\partial \pi_j}\Big(Y, X, \Delta; \pi_0 \Big)\right].\]
\end{Theorem}
The normality of Theorem \ref{thm:asynorm} on the level of the parameters $\pi$ also transfers to the quantile estimators, thanks to the multivariate delta method (see e.g. \cite{serfling_approximation_2009}). Indeed, as follows from its Theorem 3.3A, if $\sqrt{n} (\hat{\pi}_n - \pi_0) \stackrel{d}{\to} \mathcal{N}\left(0, \Sigma \right)$, then also $\sqrt{n} (g(\hat{\pi}_n) - g(\pi_0)) \stackrel{d}{\to} \mathcal{N}\left(0, D^\top\Sigma D \right)$ is implied for $D = \left(\frac{\partial g}{\partial \pi_{n,1}}, \dots, \frac{\partial g}{\partial \pi_{n,k}}\right)^\top$, whenever $g(\pi) = g(\pi_{n,1}, \dots, \pi_{n,k})$ has a nonzero differential at $\pi_0$. A sufficient condition for this is continuity of all first partial derivatives $\frac{\partial g}{\partial \pi_{n,j}}$ (for $j = 1, \dots, k$) and at least one of them being nonzero \cite[Remark 3.3B (i)]{serfling_approximation_2009}. To apply this result to our case, for any quantile level $p \in (0,1)$ and covariate value $X = x$, define the functions $h_p(\tilde{\phi}, \phi, \lambda) = Q_{\eal(\tilde{\phi}, \phi, \lambda)}(p)$ and
\begin{equation}\label{eq:deltamethod}
	g_{p,x}(\pi) = x^\top \beta + \sigma(x; \gamma) Q_{\eal(\tilde{\phi}, \phi, \lambda)}(p) = x^\top \beta + \sigma(x; \gamma) h_p(\tilde{\phi}, \phi, \lambda).
\end{equation}
Then $Q_{T|X; \theta_T} (p|x) = g_{p,x}(\pi)$ in view of relation \eqref{eq:qushift} between the quantiles of $T|X$ and the parameters of its distribution, and hence Theorem \ref{thm:asynorm} implies asymptotic normality of the quantiles, too, whenever the conditions for the delta method are satisfied.

\begin{Corollary}\label{cor:deltamethod}
	Suppose that $p \in (0,1)$ is such that the differential of $h_p(\tilde{\phi}, \phi, \lambda)$ exists. Moreover assume that for a specified value $X = x$, the resulting differential $D = \left(\frac{\partial g_{p,x}}{\partial \beta_0}, \dots, \frac{\partial g_{p,x}}{\partial \gamma_p}\right)^\top$ is nonzero at $\pi_0$. Then, under the assumptions of Theorem \ref{thm:asynorm} together with \ref{ass:norm:deltamethod}, also
	\[\sqrt{n} \Big(\hat{Q}_{n}(p|x) - Q_{T|X; \theta_T}(p|x)\Big) \stackrel{d}{\to} \mathcal{N}\left(0, D^\top (\Gamma^\top \Gamma)^{-1} \Gamma^\top V \Gamma (\Gamma^\top \Gamma)^{-1} D \right),\]
	where the estimator $\hat{Q}_n$ is obtained by replacing every parameter on the right-hand side in \eqref{eq:qushift} by the corresponding component of the estimator $\hat{\pi}_n$.
\end{Corollary}

The derivation of Corollary \ref{cor:deltamethod} as well as a remark on when its conditions are satisfied can be found in Supplement S2.4.1. An easy sufficient condition for the nonvanishing is the following:

\begin{Remark}
	All components of $\beta$ occur in the first term of \eqref{eq:deltamethod} only, such that those partial derivatives equal the corresponding components of the covariate $\tilde{x}$. A sufficient condition for the nonvanishing part (assuming that all partial derivatives -- also with respect to other parameter components than $\beta$ -- are continuous at $\pi_0$), is hence that at least one component of $\tilde{x}$ be nonzero.
\end{Remark}

Finally, we point out that the asymptotic results in this section have been derived under the assumption of fixed degrees. As will be discussed in the next section, we will relax this in the numerical implementation for flexibility. For the variance estimation in particular, implementing these asymptotic formulae in the numerical procedure would lead to underestimation due to ignorance of the variability introduced by the degree selection. A more accurate variance estimate can therefore be obtained by bootstrapping, where the optimal Laguerre degrees are determined in each bootstrap replication.

\section{Numerical implementation}\label{sec:parest}
To assess the finite sample behaviour of the estimators in a simulation setup (Section \ref{sec:sim}) and illustrate the model performance by means of a real data example (Section \ref{sec:realdata}), we still need a translation of our theoretical model to a numerical parameter estimation procedure. The completely parametric model allows for maximum likelihood optimisation, though some model aspects induce a quite challenging numerical estimation procedure. Continuity constraints are to be imposed on the parameters, and the resulting likelihood is still nondifferentiable as well as nonconvex. Also the selection of the degrees $(\tilde{m}, m)$ of the Laguerre polynomials in the model -- fixed in the theory so far -- needs to be dealt with. Finally, due to the nature of our model, the number of parameters is large, thus further complicating the estimation. A multistep algorithm tailored to the problem characteristics was conceived in \texttt{R}. The key points are explained here; more details on the initial value generation can be found in Supplement S3.
\bigbreak
\textbf{Laguerre polynomial degree: AIC selection.} The choice of the degree of the Laguerre polynomials is made based on the data at hand, such that the degrees are sufficiently high to ensure enough flexibility to capture the data trends, while on the other hand restricting the amount of parameters for a feasible estimation procedure. The Akaike information criterion (AIC) is determined for models fitted over a grid of degree pairs $(\tilde{m}, m)$ (not necessarily equal); the default upper bound is set to $4$ for both entries. Those Laguerre degrees are selected for which the AIC criterion is minimised. Simulation studies showed that the maximal degrees are almost never taken on, indicating that this default upper bound is not too strict.
\bigbreak
\textbf{Nonconvexity: multiple initial values.} Since the likelihood to be optimised is highly nonconvex, multiple starting values are chosen in order to avoid getting stuck in local maxima. These are chosen in a mostly data-driven way for the first step; in the subsequent steps, they are generated by slightly perturbing the resulting parameters of the previous step (cf. Supplement S3). The number of initial values is mostly determined by the computational cost per algorithm step.
\bigbreak
\textbf{Nondifferentiability and constraints: hybrid optimisation.} The two-sided definition \eqref{eq:dens:eal} implies that the EAL density is, in general, not even continuous at the origin, hence the likelihood is certainly nondifferentiable on the regression curve, containing the singularities $y_i = x_i^\top \beta$ (see also \eqref{eq:Tdistr}). In a similar derivative-free, nonconvex context, \cite{ewnetu_flexible_2023} opted for the Nelder-Mead (NM) algorithm because of its good performance despite the circumstances. Next, even if a not necessarily continuous likelihood behaviour does not immediately yield theoretical problems for the quantile estimation, it is more natural to impose continuity on the used EAL density, especially when taking on the distributional regression point of view (cf. Remark \ref{rem:model}\ref{rem:model:distr.reg}). Using the expression for the polynomial coefficients of the Laguerre polynomials given in (\ref{eq:lag:coeff}), continuity at the origin comes down to the following condition on $\tilde{\phi}$ and $\phi$:
\begin{equation}\label{eq:phi:ctu:cond}
	\frac{(1 + \tilde{\phi}_1 + \dots + \tilde{\phi}_{\tilde{m}})^2}{1 + \tilde{\phi}_1^2 + \dots + \tilde{\phi}_{\tilde{m}}^2} = \frac{(1 + \phi_1 + \dots + \phi_m)^2}{1 + \phi_1^2 + \dots + \phi_m^2}.
\end{equation}
Observe that symmetry ($\tilde{\phi} = \phi$) is a sufficient, but not a necessary condition. As the NM algorithm is unable to handle constrained optimisation with constraints of the form \eqref{eq:phi:ctu:cond}, we constructed a two-step optimisation method \texttt{NMCob} that concatenates the stability of the NM algorithm (\texttt{R} package \texttt{optim}) with the ability of \texttt{COBYLA} (\texttt{nloptr} package) to deal with constraints as in \eqref{eq:phi:ctu:cond}. More weight (400 iterations) is given to the NM part -- better in finding its way towards the optimum, even when still far away-- while the final (100) iterations serve only to impose the continuity, although the remaining parameters of $T$ are naturally allowed to vary with the coefficients $\tilde{\phi}$ and $\phi$.
\bigbreak
\textbf{Summary: multistep approach.} The algorithm consists of three steps. In the basis step, optimisation over all parameters simultaneously is done putting the Laguerre degrees to zero. This reduces the number of parameters and results in more stable estimates for the remaining parameters (in particular for the regression coefficient $\beta$). Since constraint \eqref{eq:phi:ctu:cond} becomes void under the absence of the Laguerre polynomials, only NM (500 iterations) rather than NMCob suffices. Next, in the intermediate step, our two-step optimisation method NMCob is applied over a grid of degrees to select the optimal Laguerre degrees. In order to somewhat reduce the large computational cost, optimisation is done over the parameters for $T$ only, while fixing those for the censoring $C$ and the copula to the results of the basis step. In the final step, another simultaneous optimisation over all parameters takes place, this time with the polynomial degrees fixed to the selected optimal values, and hence again using NMCob rather than NM.

\section{Simulations}\label{sec:sim}
Three major simulation scenarios are considered, each of them with a number of subscenarios. Each time, 200 samples are considered, all of size 500 unless mentioned otherwise. In Scenario 1, the most extensive one, we study the impact of some characteristics of both the generated data and the model fitting procedure. It moreover illustrates that applying the model when the independence of the data at hand is uncertain, can be considered quite safe: while still performing well in the independent case, it greatly improves the modelling results under actual dependent censoring. This first scenario is treated in the main text below with a focus on the quantile performance; the parameter performance can be found in Supplement S4.1. The Supplementary material also contains a detailed description of Scenario 2 (Supplement S4.2) and Scenario 3 (Supplement S4.3) and a discussion of their results; a brief conclusion is included below (Section \ref{sec:sim:scen23}).

In each of the scenarios, the censoring is correctly specified as a normal distribution. This is to avoid an endless amount of different scenarios next to the options for both the copula and the survival time distribution that are already there (including misspecification in both these components simultaneously). Moreover, \cite{deresa_copula_2024} performed simulation studies in a similar, dependent censoring setup and concluded in their Section 6 that misspecification of the censoring (or even both the censoring and the copula) doesn't negatively influence the parameter estimates for the survival distribution too much.

\subsection{Scenario 1}
In this first scenario, both the (logarithmic) survival time $T$ and the (logarithmic) censoring time $C$ are generated from a normal distribution:
\begin{equation}\label{eq:sim:scen1}
	\begin{aligned}
		T &= \beta_0 + \beta_1\tilde{X} + \exp(\gamma_0 + \gamma_1\tilde{X})\epsilon_T, &\qquad \epsilon_T \sim \mathcal{N}(0,1), \\
		C &= \alpha_0 + \alpha_1\tilde{X} + \sigma_C \epsilon_C, &\qquad \epsilon_C \sim \mathcal{N}(0,1).
	\end{aligned}
\end{equation}
The covariate $\tilde{X}$ is uniformly distributed over $[0,4]$. For the dependence structure, a Frank copula is used for both the data generation and the modelling. The copula parameter for the generation is chosen corresponding to a Kendall's tau of $\tau_K = 0.5$. The remaining parameters are chosen in such a way that the proportion of uncensored observations equals approximately 50-55\%: we put $(\beta_0, \beta_1) = (2.8, 0.6)$, $(\gamma_0, \gamma_1) = (-1.5, 0.45)$, $(\alpha_0, \alpha_1) = (3.15, 0.45)$ and $\sigma_C = 0.8$. Only in the homoscedastic subscenario, the parameter $\lambda$ is fixed during the modelling to ensure identifiability (cf. Section \ref{sec:model:multiqu:parastatus:lam}). In this case, $(\gamma_0, \gamma_1)$ is replaced by one parameter $\gamma_T = -1.7$, such that $\sigma_T = \exp(\gamma_T) \approx 0.18$ (or, equivalently, put $\sigma_T(\tilde{X}, \gamma_0, \gamma_1) = \exp(\gamma_0 + \gamma_1\tilde{X})$ with, in the homoscedastic case, $\gamma_0 = -1.7$ and $\gamma_1 = 0$). The other parameters for $T$ are never changed, and for $C$ only the intercept $\alpha_0$ is subject to some variation in order to obtain the desired censoring rate.

An overview of all subscenarios is given in Table \ref{tab:scen1:subscen} below, where the top rows contain the heteroscedastic basis scenario, and its homoscedastic counterpart with fixed value $\lambda = 0.5$. Note that for the homoscedastic case we switched to the Clayton family only to illustrate that this family also works; we could just as well have preserved the Frank copula here. (Comparing the results of the hetero- and homoscedastic scenario would not be possible anyway, even if we would have used the same copula families, since the model fitted is really different.) First, we study how important the censoring proportion is, how the results vary under different sample sizes, and whether the strength of the dependence during the data generation significantly impacts the model performance. Next, we impose a restriction on the parameter of the Frank copula in the model fitting phase such that only positive dependence values (corresponding $\tau_K \in [0,1]$) can be obtained. This enables seeing whether knowing the dependence \enquote*{direction} beforehand makes a big difference. Afterwards, we assess the influence of the fixed value for $\lambda$ used during the modelling in a homoscedastic setup. In the following series of settings, the data are either generated or modelled using an independence copula (or both). This allows to see the bias one would introduce when the underlying dependence would be ignored and, on the other hand, it illustrates that using the dependent copula model is not harmful even when there is actually no dependence (and it can be used also in the case where independence is correctly assumed). Finally, we study the impact of copula misspecification.

\begin{table}[ht]
	\centering
	\caption{Subscenarios considered within Scenario 1. All \enquote*{-} entries are equal to the \enquote*{BasisHet} scenario, only differences are explicitly specified. The slash in the column for $\lambda$ (in all heteroscedastic scenarios) means that this parameter is not subject to any restrictions, whereas in the homoscedastic scenario it remains fixed during the modelling. Abbreviations used: \enquote*{UnCens} = approximate proportion of uncensored observations, \enquote*{gen.} = during the data generation process, \enquote*{fit} = for the model fitting, \enquote*{Hom} = homoscedastic, \enquote*{Het} = heteroscedastic, \enquote*{indep} = independence (copula), \enquote*{MSCop} = misspecified copula family, \enquote*{FrankPos} = Frank copula, with additional restriction for positive dependence only. ($\dagger$) Recall that in the homoscedastic scenarios, $\gamma_0 = -1.7$ and $\gamma_1 = 0$, rather than the basis value for $\gamma$. ($\ddagger$) In order to obtain these uncensoring rates, the basis value of $\alpha_0 = 3.15$ is changed to 3.5 or 2.85 for less or more censoring, respectively.}
	\label{tab:scen1:subscen}
	\begin{tabular}{lcccc@{\hskip 25 pt}cc}
		\toprule
		\textbf{Scenario} & \multicolumn{4}{c}{\textbf{Data generation}} & \multicolumn{2}{c}{\textbf{Modelling}} \\
		& Sample size & UnCens & copula (gen.) & $\tau_K$ & copula (fit) & $\lambda$ \\
		\midrule
		BasisHet & 500 & 0.54 & Frank & 0.5 & Frank & / \\
		BasisHom$^\dagger$ & - & 0.51 & Clayton & - & Clayton & 0.5 \\
		\midrule
		LessCens$^\ddagger$ & - & 0.74 & - & - & - & / \\
		MoreCens$^\ddagger$ & - & 0.35 & - & - & -  & / \\
		\midrule
		Size S & 250 & - & - & - & - & / \\
		Size L & 1000 & - & - & - & - & / \\
		Size XL & 2000 & - & - & - & - & / \\
		\midrule
		LessDep & - & - & - & 0.25 & - & / \\
		MoreDep & - & - & - & 0.75 & - & / \\
		\midrule
		FitPos & - & - & - & - & FrankPos & / \\
		\midrule
		Hom0.3$^\dagger$ & - & 0.51 & Clayton & - & Clayton & 0.3 \\
		Hom0.7$^\dagger$ & - & 0.51 & Clayton & - & Clayton & 0.7 \\
		\midrule
		FitIndep & - & - & Frank & - & indep  & / \\
		GenIndep & - & 0.53 & indep & 0 & Frank  & / \\
		AllIndep  & - & 0.53 & indep & 0 & indep  & / \\
		FitIndepHom$^\dagger$ & - & 0.51 & Clayton & - & indep  & 0.5 \\
		\midrule
		MSCopHet & - & - & Frank & - & Gumbel & / \\
		MSCopHom$^\dagger$ & - & 0.51 & Clayton & - & Gumbel & 0.5 \\
		\bottomrule
	\end{tabular}
\end{table}

The performance of the model in all these subscenarios is assessed both on the level of the parameters and on the level of the quantiles. In Supplement S4.1, we look at how well the parameters themselves are estimated, which is already fairly indicative of the performance on the level of the quantiles as well. Nonetheless, as discussed there, comparison to the true parameter from the data generation process is not always meaningful or even impossible for some of the parameters, especially those of the EAL distribution. Since also these parameters exert an important influence on the overall fit of the model, it is informative to have a look at the quantile predictions for some specified covariate values and quantile levels.

\subsubsection{Quantile performance}\label{sec:sim:scen1:quEst}
We consider the three quartile levels $p = 0.25$, $p = 0.50$ and $p = 0.75$, each time for three values of the uniformly distributed covariate $\tilde{X} \sim \mathcal{U}[0,4]$, namely $\tilde{x} = 1$, $\tilde{x} = 2$ and $\tilde{x} = 3$. The true values are reported together with the average estimate over all 200 samples, the empirical variance and finally the relative bias in Tables \ref{tab:scen1:qu:res:basis} through \ref{tab:scen1:qu:res:MSCop}. (The empirical variance is multiplied by ten in the tables in order to include an additional digit, since typically the first digits are zero. We report the empirical instead of bootstrap-based variance estimate to reduce the already high computational cost; why these provide a more realistic variance estimate than the asymptotic formulae of Section \ref{sec:asy} has been discussed there.). The results for the basis scenarios are presented in Table \ref{tab:scen1:qu:res:basis} below. All remaining tables contain slight variations and are to be compared to the former in order to study the impact of some of the model aspects.

\begin{table}[ht]
	\centering
	\caption{Quantile performance for Scenario 1 (basis). Reported are the true values, average estimate (avg.), empirical variance multiplied by ten (eVar$\times 10$) and relative bias (rBias).}
	\label{tab:scen1:qu:res:basis}
	\begin{adjustbox}{width = \textwidth}
		\begin{tabular}{llccc@{\hskip 25 pt}ccc@{\hskip 25 pt}ccc}
			\toprule
			\textbf{Scenario} & & \multicolumn{3}{c}{$p = 0.25$} & \multicolumn{3}{c}{$p = 0.50$} & \multicolumn{3}{c}{$p = 0.75$} \\
			& & $\tilde{x} = 1$ & $\tilde{x} = 2$ & $\tilde{x} = 3$ & $\tilde{x} = 1$ & $\tilde{x} = 2$ & $\tilde{x} = 3$ & $\tilde{x} = 1$ & $\tilde{x} = 2$ & $\tilde{x} = 3$ \\
			\midrule
			\multirow{4}{*}{BasisHet}
			& true & 3.164 & 3.630 & 4.020 & 3.400 & 4.000 & 4.600 & 3.636 & 4.370 & 5.181 \\ 
			& avg. & 3.206 & 3.688 & 4.088 & 3.435 & 4.055 & 4.680 & 3.659 & 4.416 & 5.260 \\ 
			& eVar$\times 10$ & 0.036 & 0.069 & 0.108 & 0.032 & 0.073 & 0.142 & 0.033 & 0.090 & 0.238 \\ 
			& rBias & 0.013 & 0.016 & 0.017 & 0.010 & 0.014 & 0.017 & 0.006 & 0.010 & 0.015 \\ 
			\midrule
			\multirow{4}{*}{BasisHom}
			& true & 3.277 & 3.877 & 4.477 & 3.400 & 4.000 & 4.600 & 3.523 & 4.123 & 4.723 \\ 
			& avg. & 3.295 & 3.898 & 4.501 & 3.410 & 4.014 & 4.617 & 3.526 & 4.129 & 4.733 \\ 
			& eVar$\times 10$ & 0.011 & 0.011 & 0.014 & 0.006 & 0.006 & 0.008 & 0.004 & 0.004 & 0.006 \\ 
			& rBias & 0.006 & 0.006 & 0.006 & 0.003 & 0.003 & 0.004 & 0.001 & 0.002 & 0.002 \\ 
			\bottomrule
		\end{tabular}
	\end{adjustbox}
\end{table}

\textbf{General discussion.} Let us first point out that the performance is good in both the heteroscedastic (\enq{BasisHet}) as well as the homoscedastic (\enq{BasisHom}) basis case (Table \ref{tab:scen1:qu:res:basis}). This is reassuring: it shows that our model, fixing $\lambda$ depending on the situation, works in both cases, and we confirm that indeed $\lambda$ does not have to be fixed in the heteroscedastic case in order to obtain identifiability.

We proceed by making a few straightforward observations that are in line with what one would expect. To begin with, the quantile performance drastically decreases when the amount of censoring increases. The \enquote*{LessCens} case is better than \enq{BasisHet} both in terms of the bias and in terms of the variance, while \enq{MoreCens} is significantly worse (Table \ref{tab:scen1:qu:res:censprop}). Conversely, the performance increases with the sample size (Table \ref{tab:scen1:qu:res:size}). At some point, however, increasing the sample size even more, barely influences the performance in terms of the bias, as can be seen by comparing the almost identical relative bias of L and XL sample sizes. Nevertheless, the variance still decreases with increasing sample size, as it should. Whereas the amount of censoring does influence the performance, the amount of dependence apparently does so only to a limited extent, mainly in terms of a reduced variance for the \enq{MoreDep} scenario (Table \ref{tab:scen1:qu:res:DepStrength}). As for knowing the dependence direction (Table \ref{tab:scen1:qu:res:DepDir}), surprisingly this does not really improve the results. A possible explanation for this could be that the actual dependence ($\tau_K = 0.5$) is sufficiently positive, such that also the algorithm without restrictions immediately starts exploring the right area and very little iterations are lost to copula parameter regions corresponding to negative dependence.

\begin{table}[ht]
	\centering
	\caption{Quantile performance for Scenario 1 (censoring proportion). Reported are the true values, average estimate (avg.), empirical variance multiplied by ten (eVar$\times 10$) and relative bias (rBias).}
	\label{tab:scen1:qu:res:censprop}
	\begin{adjustbox}{width = \textwidth}
		\begin{tabular}{llccc@{\hskip 25 pt}ccc@{\hskip 25 pt}ccc}
			\toprule
			\textbf{Scenario} & & \multicolumn{3}{c}{$p = 0.25$} & \multicolumn{3}{c}{$p = 0.50$} & \multicolumn{3}{c}{$p = 0.75$} \\
			& & $\tilde{x} = 1$ & $\tilde{x} = 2$ & $\tilde{x} = 3$ & $\tilde{x} = 1$ & $\tilde{x} = 2$ & $\tilde{x} = 3$ & $\tilde{x} = 1$ & $\tilde{x} = 2$ & $\tilde{x} = 3$ \\
			& \textbf{true} & 3.164 & 3.630 & 4.020 & 3.400 & 4.000 & 4.600 & 3.636 & 4.370 & 5.181 \\ 
			\midrule
			\multirow{3}{*}{LessCens}
			& avg. & 3.177 & 3.652 & 4.052 & 3.409 & 4.019 & 4.632 & 3.638 & 4.380 & 5.204 \\ 
			& eVar$\times 10$ & 0.016 & 0.025 & 0.051 & 0.011 & 0.020 & 0.051 & 0.012 & 0.027 & 0.081 \\ 
			& rBias & 0.004 & 0.006 & 0.008 & 0.003 & 0.005 & 0.007 & 0.001 & 0.002 & 0.005 \\
			\midrule
			\multirow{3}{*}{MoreCens}
			& avg. & 3.275 & 3.791 & 4.218 & 3.495 & 4.151 & 4.812 & 3.702 & 4.491 & 5.374 \\ 
			& eVar$\times 10$ & 0.124 & 0.259 & 0.371 & 0.089 & 0.211 & 0.341 & 0.071 & 0.194 & 0.412 \\ 
			& rBias & 0.035 & 0.044 & 0.049 & 0.028 & 0.038 & 0.046 & 0.018 & 0.028 & 0.037 \\ 
			\bottomrule
		\end{tabular}
	\end{adjustbox}
\end{table}

\begin{table}[ht]
	\centering
	\caption{Quantile performance for Scenario 1 (sample size). Reported are the true values, average estimate (avg.), empirical variance multiplied by ten (eVar$\times 10$) and relative bias (rBias).}
	\label{tab:scen1:qu:res:size}
	\begin{adjustbox}{width = \textwidth}
		\begin{tabular}{llccc@{\hskip 25 pt}ccc@{\hskip 25 pt}ccc}
			\toprule
			\textbf{Scenario} & & \multicolumn{3}{c}{$p = 0.25$} & \multicolumn{3}{c}{$p = 0.50$} & \multicolumn{3}{c}{$p = 0.75$} \\
			& & $\tilde{x} = 1$ & $\tilde{x} = 2$ & $\tilde{x} = 3$ & $\tilde{x} = 1$ & $\tilde{x} = 2$ & $\tilde{x} = 3$ & $\tilde{x} = 1$ & $\tilde{x} = 2$ & $\tilde{x} = 3$ \\
			& \textbf{true} & 3.164 & 3.630 & 4.020 & 3.400 & 4.000 & 4.600 & 3.636 & 4.370 & 5.181 \\
			\midrule
			\multirow{3}{*}{Size S}
			& avg. & 3.216 & 3.710 & 4.131 & 3.439 & 4.069 & 4.712 & 3.658 & 4.422 & 5.288 \\ 
			& eVar$\times 10$ & 0.067 & 0.114 & 0.170 & 0.065 & 0.135 & 0.293 & 0.100 & 0.268 & 0.900 \\ 
			& rBias & 0.017 & 0.022 & 0.028 & 0.012 & 0.017 & 0.024 & 0.006 & 0.012 & 0.021 \\ 
			\midrule
			\multirow{3}{*}{Size L} 
			& avg. & 3.187 & 3.665 & 4.075 & 3.421 & 4.037 & 4.670 & 3.643 & 4.392 & 5.238 \\ 
			& eVar$\times 10$ & 0.023 & 0.043 & 0.067 & 0.017 & 0.043 & 0.089 & 0.012 & 0.035 & 0.100 \\ 
			& rBias & 0.007 & 0.010 & 0.014 & 0.006 & 0.009 & 0.015 & 0.002 & 0.005 & 0.011 \\ 
			\midrule
			\multirow{3}{*}{Size XL} 
			& avg. & 3.194 & 3.671 & 4.074 & 3.423 & 4.038 & 4.665 & 3.646 & 4.396 & 5.239 \\ 
			& eVar$\times 10$ & 0.012 & 0.021 & 0.032 & 0.009 & 0.022 & 0.050 & 0.007 & 0.021 & 0.061 \\ 
			& rBias & 0.009 & 0.011 & 0.014 & 0.007 & 0.010 & 0.014 & 0.003 & 0.006 & 0.011 \\ 
			\bottomrule
		\end{tabular}
	\end{adjustbox}
\end{table}
\bigbreak
\textbf{Bias trends: varying covariate values.} Observe that within a triplet of columns corresponding to the same quantile level (e.g. the first three columns for $p = 0.25$ in Table \ref{tab:scen1:qu:res:size}), the relative bias is increasing for larger covariate values and gets pretty high for $\tilde{x} = 3$ sometimes. This could be explained by the fact that also the censoring proportion increases with $\tilde{x}$. Indeed, as $\beta_1 > \alpha_1$, i.e. the slope for $T$'s regression line is larger than that for $C$, so the probability that the value of $T$ corresponding to a high value of $\tilde{X}$ lies above the censoring value $C$ (and hence is no longer observed) is increasing. Therefore, the rightmost part of the covariate distribution is subject to more heavy censoring, which in turn results in worse quantile prediction performance.

\bigbreak
\textbf{Bias trends: varying quantile levels.} A bias trend can be observed also in the quantile level of interest. Comparing relative bias values for the same covariate value (e.g. $\tilde{x} = 1$) over all three quantile levels (so columns 1, 4 and 7) in Table \ref{tab:scen1:qu:res:size}, for instance, we note that the bias decreases with $p$. This is remarkable, as in the right tail of $T$'s distribution we have very little observations, due to the censoring. Therefore, one would expect the opposite of the presently observed trend. The exact reason for this is hard to infer. It could be due to the model misspecification (on the level of $T$ and, in some cases, the copula, too), inducing a different bias behaviour than expected.

\begin{table}[ht]
	\centering
	\caption{Quantile performance for Scenario 1 (dependence strength). Reported are the true values, average estimate (avg.), empirical variance multiplied by ten (eVar$\times 10$) and relative bias (rBias).}
	\label{tab:scen1:qu:res:DepStrength}
	\begin{adjustbox}{width = \textwidth}
		\begin{tabular}{llccc@{\hskip 25 pt}ccc@{\hskip 25 pt}ccc}
			\toprule
			\textbf{Scenario} & & \multicolumn{3}{c}{$p = 0.25$} & \multicolumn{3}{c}{$p = 0.50$} & \multicolumn{3}{c}{$p = 0.75$} \\
			& & $\tilde{x} = 1$ & $\tilde{x} = 2$ & $\tilde{x} = 3$ & $\tilde{x} = 1$ & $\tilde{x} = 2$ & $\tilde{x} = 3$ & $\tilde{x} = 1$ & $\tilde{x} = 2$ & $\tilde{x} = 3$ \\
			& \textbf{true} & 3.164 & 3.630 & 4.020 & 3.400 & 4.000 & 4.600 & 3.636 & 4.370 & 5.181 \\
			\midrule
			\multirow{3}{*}{LessDep} 
			& avg. & 3.177 & 3.648 & 4.043 & 3.411 & 4.019 & 4.633 & 3.648 & 4.395 & 5.230 \\ 
			& eVar$\times 10$ & 0.032 & 0.065 & 0.106 & 0.033 & 0.085 & 0.187 & 0.038 & 0.109 & 0.309 \\ 
			& rBias & 0.004 & 0.005 & 0.006 & 0.003 & 0.005 & 0.007 & 0.003 & 0.006 & 0.010 \\ 
			\midrule
			\multirow{3}{*}{MoreDep} 
			& avg. & 3.197 & 3.667 & 4.057 & 3.411 & 4.019 & 4.638 & 3.631 & 4.381 & 5.235 \\ 
			& eVar$\times 10$ & 0.016 & 0.029 & 0.073 & 0.008 & 0.018 & 0.051 & 0.009 & 0.017 & 0.0074 \\ 
			& rBias & 0.010 & 0.010 & 0.009 & 0.003 & 0.005 & 0.008 & -0.001 & 0.003 & 0.011 \\ 
			\bottomrule
		\end{tabular}
	\end{adjustbox}
\end{table}

\begin{table}[ht]
	\centering
	\caption{Quantile performance for Scenario 1 (dependence direction). Reported are the true values, average estimate (avg.), empirical variance multiplied by ten (eVar$\times 10$) and relative bias (rBias).}
	\label{tab:scen1:qu:res:DepDir}
	\begin{adjustbox}{width = \textwidth}
		\begin{tabular}{llccc@{\hskip 25 pt}ccc@{\hskip 25 pt}ccc}
			\toprule
			\textbf{Scenario} & & \multicolumn{3}{c}{$p = 0.25$} & \multicolumn{3}{c}{$p = 0.50$} & \multicolumn{3}{c}{$p = 0.75$} \\
			& & $\tilde{x} = 1$ & $\tilde{x} = 2$ & $\tilde{x} = 3$ & $\tilde{x} = 1$ & $\tilde{x} = 2$ & $\tilde{x} = 3$ & $\tilde{x} = 1$ & $\tilde{x} = 2$ & $\tilde{x} = 3$ \\
			\midrule
			\multirow{4}{*}{FitPos}
			& true & 3.164 & 3.630 & 4.020 & 3.400 & 4.000 & 4.600 & 3.636 & 4.370 & 5.181 \\ 
			& avg. & 3.208 & 3.695 & 4.098 & 3.437 & 4.063 & 4.692 & 3.662 & 4.426 & 5.279 \\ 
			& eVar$\times 10$ & 0.030 & 0.061 & 0.102 & 0.029 & 0.072 & 0.138 & 0.026 & 0.070 & 0.173 \\ 
			& rBias & 0.014 & 0.018 & 0.020 & 0.011 & 0.016 & 0.020 & 0.007 & 0.013 & 0.019 \\
			\bottomrule
		\end{tabular}
	\end{adjustbox}
\end{table}

\bigbreak
\textbf{The choice of $\lambda$.} Although the homoscedastic model is, in theory, equally identifiable for any value of $\lambda$ (as discussed in Section \ref{sec:model:multiqu:parastatus:lam}), we note that in practice some values seem to be better than others for converging to the true quantiles. In this Scenario 1, $\lambda = 0.5$ is clearly the better value (cf. Tables \ref{tab:scen1:qu:res:basis} and \ref{tab:scen1:qu:res:Hom}), corresponding to the smallest bias for each of the estimated quantiles; on the level of the variance, there seems to be no significant difference. We hypothesise that values of $\lambda$ for which the basis shape of the AL distribution is closer to $T$’s true distribution, lead to overall better results. This could be due to the constraint on the Laguerre polynomials: in practice, we always impose an upper bound on the degree. Consequently, $\lambda$-values for which the corresponding basis AL distribution (i.e. with Laguerre polynomial degrees equal to zero) is closer to the shape of the true distribution, are more likely to get close to this true distribution even within the limited amount of Laguerre degree possibilities. In this specific scenario, where the true distribution for $T$ is symmetric, $\lambda = 0.5$ is not unexpectedly the best performing $\lambda$-value. However, $\lambda = 0.3$ performs consistently worse than $\lambda = 0.7$; we do not have a clear explanation for this (see also Supplement S4 for a more in-depth investigation). We therefore refrain from formulating general guidelines on how to choose $\lambda$ in practice. We advise the user of our homoscedastic model to just try out multiple values and then perform model selection. This is also illustrated in the real data application of Section \ref{sec:realdata}.

\begin{table}[ht]
	\centering
	\caption{Quantile performance for Scenario 1 (homoscedastic, $\lambda$ influence). Reported are the true values, average estimate (avg.), empirical variance multiplied by ten (eVar$\times 10$) and relative bias (rBias).}
	\label{tab:scen1:qu:res:Hom}
	\begin{adjustbox}{width = \textwidth}
		\begin{tabular}{llccc@{\hskip 25 pt}ccc@{\hskip 25 pt}ccc}
			\toprule
			\multicolumn{2}{l}{\textbf{Scenario}} & \multicolumn{3}{c}{$p = 0.25$} & \multicolumn{3}{c}{$p = 0.50$} & \multicolumn{3}{c}{$p = 0.75$} \\
			& & $\tilde{x} = 1$ & $\tilde{x} = 2$ & $\tilde{x} = 3$ & $\tilde{x} = 1$ & $\tilde{x} = 2$ & $\tilde{x} = 3$ & $\tilde{x} = 1$ & $\tilde{x} = 2$ & $\tilde{x} = 3$ \\
			& \textbf{true} & 3.277 & 3.877 & 4.477 & 3.400 & 4.000 & 4.600 & 3.523 & 4.123 & 4.723 \\
			\midrule
			\multirow{3}{*}{Hom0.3}  
			& avg. & 3.355 & 3.964 & 4.574 & 3.452 & 4.061 & 4.670 & 3.559 & 4.169 & 4.778 \\ 
			& eVar$\times 10$ & 0.008 & 0.010 & 0.014 & 0.006 & 0.007 & 0.010 & 0.005 & 0.006 & 0.009 \\ 
			& rBias & 0.024 & 0.023 & 0.022 & 0.015 & 0.015 & 0.015 & 0.010 & 0.011 & 0.012 \\ 
			\midrule
			\multirow{3}{*}{Hom0.7} 
			& avg. & 3.229 & 3.825 & 4.422 & 3.384 & 3.981 & 4.578 & 3.520 & 4.116 & 4.713 \\ 
			& eVar$\times 10$ & 0.018 & 0.020 & 0.025 & 0.005 & 0.006 & 0.009 & 0.005 & 0.005 & 0.008 \\ 
			& rBias & -0.015 & -0.013 & -0.012 & -0.005 & -0.005 & -0.005 & -0.001 & -0.002 & -0.002 \\
			\bottomrule
		\end{tabular}
	\end{adjustbox}
\end{table}

\bigbreak
\textbf{Dependence vs. independence, and copula misspecification}. This model being designed to tackle dependent censoring, it is interesting to see how it actually behaves under different dependence specifications. And even if the \emph{strength} of the dependence does not seem to be of major impact to the model performance (cf. \emph{supra}), acknowledging its \emph{presence} is very important. This is nicely illustrated in Table \ref{tab:scen1:qu:res:Indep}. If no dependence is present and estimation is done under the (true) independence copula (\enquote*{AllIndep} scenario), the model performs extremely well. But it gets more interesting when either the true copula or the fitted copula belongs to a family that actually models a dependence structure. The \enq{GenIndep} scenario, where the data are actually subject to \enquote*{standard} independent censoring, poses no problem. Applying our method yields only slightly increased bias values compared to \enq{AllIndep}, though it also comes with a moderate variance increase. Next, comparing the last row of the \enquote*{BasisHet} scenario (Table \ref{tab:scen1:qu:res:basis}) with that of \enquote*{FitIndep} (Table \ref{tab:scen1:qu:res:Indep}), in each column the relative bias in the latter case (i.e. when the dependence is ignored) amounts to at least 3.5 to 4 times the corresponding (relative) bias in the basis scenario where a Frank copula is used for the modelling, and even more than 5 times the value for $p = 0.75$, $\tilde{x} = 1$. This effect is even more pronounced in the homoscedastic case: the relative bias for \enquote*{FitIndepHom} is generally 5 to 6 times the corresponding one in \enquote*{BasisHom}, up to more than 8 times for all covariate values in the $p = 0.75$ column. Moreover, even under misspecification of the dependence structure (\enquote*{MSCopHet} scenario, Table \ref{tab:scen1:qu:res:MSCop}), the quantile bias is still less than when fitting under independence (\enquote*{FitIndep}) without a single exception, although it comes at the cost of increased variability. The same holds true when comparing \enquote*{MSCopHom} with \enquote*{FitIndepHom} (Table \ref{tab:scen1:qu:res:MSCop}). We conclude that applying our model that allows for dependence is never very harmful, and, in the case of actual dependence, greatly performance improving compared to completely ignoring the dependence, even in the case of copula misspecification.

\begin{table}[ht]
	\centering
	\caption{Quantile performance for Scenario 1 (independence). Reported are the true values, average estimate (avg.), empirical variance multiplied by ten (eVar$\times 10$) and relative bias (rBias). Only in the homoscedastic scenario, the true quantile values differ from those reported in the table header.}
	\label{tab:scen1:qu:res:Indep}
	\begin{adjustbox}{width = \textwidth}
		\begin{tabular}{llccc@{\hskip 25 pt}ccc@{\hskip 25 pt}ccc}
			\toprule
			\textbf{Scenario} & & \multicolumn{3}{c}{$p = 0.25$} & \multicolumn{3}{c}{$p = 0.50$} & \multicolumn{3}{c}{$p = 0.75$} \\
			& & $\tilde{x} = 1$ & $\tilde{x} = 2$ & $\tilde{x} = 3$ & $\tilde{x} = 1$ & $\tilde{x} = 2$ & $\tilde{x} = 3$ & $\tilde{x} = 1$ & $\tilde{x} = 2$ & $\tilde{x} = 3$ \\
			& \textbf{true} & 3.164 & 3.630 & 4.020 & 3.400 & 4.000 & 4.600 & 3.636 & 4.370 & 5.181 \\ 
			\midrule
			\multirow{3}{*}{FitIndep}
			& avg. & 3.311 & 3.845 & 4.265 & 3.539 & 4.227 & 4.907 & 3.757 & 4.594 & 5.524 \\ 
			& eVar$\times 10$ & 0.015 & 0.033 & 0.073 & 0.011 & 0.028 & 0.067 & 0.012 & 0.033 & 0.123 \\ 
			& rBias & 0.047 & 0.059 & 0.061 & 0.041 & 0.057 & 0.067 & 0.033 & 0.051 & 0.066 \\ 
			\midrule
			\multirow{3}{*}{GenIndep} 
			& avg. & 3.154 & 3.623 & 4.032 & 3.394 & 4.001 & 4.630 & 3.645 & 4.399 & 5.261 \\ 
			& eVar$\times 10$ & 0.027 & 0.057 & 0.114 & 0.049 & 0.127 & 0.382 & 0.137 & 0.420 & 1.471 \\ 
			& rBias & -0.003 & -0.002 & 0.003 & -0.002 & 0.000 & 0.007 & 0.003 & 0.007 & 0.016 \\ 
			\midrule
			\multirow{3}{*}{AllIndep} 
			& avg. & 3.169 & 3.635 & 4.028 & 3.399 & 3.995 & 4.593 & 3.631 & 4.360 & 5.168 \\ 
			& eVar$\times 10$ & 0.012 & 0.026 & 0.058 & 0.008 & 0.021 & 0.057 & 0.011 & 0.029 & 0.106 \\ 
			& rBias & 0.002 & 0.001 & 0.002 & -0.000 & -0.001 & -0.002 & -0.001 & -0.002 & -0.002 \\
			\midrule
			\multirow{4}{40pt}{FitIndep\-Hom}
			& true & 3.277 & 3.877 & 4.477 & 3.400 & 4.000 & 4.600 & 3.523 & 4.123 & 4.723 \\ 
			& avg. & 3.380 & 3.993 & 4.606 & 3.471 & 4.083 & 4.696 & 3.573 & 4.186 & 4.799 \\ 
			& eVar$\times 10$ & 0.002 & 0.001 & 0.003 & 0.003 & 0.002 & 0.003 & 0.003 & 0.003 & 0.004 \\ 
			& rBias & 0.032 & 0.030 & 0.029 & 0.021 & 0.021 & 0.021 & 0.014 & 0.015 & 0.016 \\ 
			\bottomrule
		\end{tabular}
	\end{adjustbox}
\end{table}

\begin{table}[ht]
	\centering
	\caption{Quantile performance for Scenario 1 (copula misspecification). Reported are the true values, average estimate (avg.), empirical variance multiplied by ten (eVar$\times 10$) and relative bias (rBias).}
	\label{tab:scen1:qu:res:MSCop}
	\begin{adjustbox}{width = \textwidth}
		\begin{tabular}{llccc@{\hskip 25 pt}ccc@{\hskip 25 pt}ccc}
			\toprule
			\textbf{Scenario} & & \multicolumn{3}{c}{$p = 0.25$} & \multicolumn{3}{c}{$p = 0.50$} & \multicolumn{3}{c}{$p = 0.75$} \\
			& & $\tilde{x} = 1$ & $\tilde{x} = 2$ & $\tilde{x} = 3$ & $\tilde{x} = 1$ & $\tilde{x} = 2$ & $\tilde{x} = 3$ & $\tilde{x} = 1$ & $\tilde{x} = 2$ & $\tilde{x} = 3$ \\
			\midrule
			\multirow{4}{*}{MSCopHet} 
			& true & 3.164 & 3.630 & 4.020 & 3.400 & 4.000 & 4.600 & 3.636 & 4.370 & 5.181 \\ 
			& avg. & 3.255 & 3.763 & 4.175 & 3.478 & 4.126 & 4.766 & 3.680 & 4.454 & 5.301 \\ 
			& eVar$\times 10$ & 0.037 & 0.079 & 0.134 & 0.035 & 0.090 & 0.188 & 0.034 & 0.103 & 0.294 \\ 
			& rBias & 0.029 & 0.037 & 0.039 & 0.023 & 0.032 & 0.036 & 0.012 & 0.019 & 0.023 \\ 
			\midrule
			\multirow{4}{*}{MSCopHom}
			& true & 3.277 & 3.877 & 4.477 & 3.400 & 4.000 & 4.600 & 3.523 & 4.123 & 4.723 \\ 
			& avg. & 3.354 & 3.960 & 4.565 & 3.442 & 4.048 & 4.654 & 3.540 & 4.146 & 4.752 \\ 
			& eVar$\times 10$ & 0.005 & 0.005 & 0.009 & 0.005 & 0.006 & 0.009 & 0.005 & 0.005 & 0.008 \\  
			& rBias & 0.024 & 0.021 & 0.020 & 0.013 & 0.012 & 0.012 & 0.005 & 0.006 & 0.006 \\
			\bottomrule
		\end{tabular}
	\end{adjustbox}
\end{table}

\subsection{Scenario 2 \& 3}\label{sec:sim:scen23}
A more in-depth discussion for the remaining two simulation scenarios can be found in the Supplementary material (Supplement S4.2 and S4.3). The purpose of Scenario 2 is twofold. On the one hand we study the model performance when all components in the model are correctly specified in both a heteroscedastic and a homoscedastic setting, and under another covariate distribution than in Scenario 1. On the other hand, fixing $\lambda$ to different values in the homoscedastic case, we continue our investigation on how this choice influences the simulation results. In the final, fully homoscedastic Scenario 3, we go into even more depth on the influence of the fixed value for $\lambda$. Rather than the admittedly somewhat artificial setup of Scenario 2 with an EAL distribution for $T$, the more natural, yet still asymmetric, Weibull case is considered.

As already hinted at before, the behaviour in terms of the fixed value for $\lambda$ remains mysterious; we can once again only conclude that upon application, the homoscedastic model should be fitted for several $\lambda$-values and afterwards model selection is to be performed. Another surprising outcome is that the results under the fully correctly specified model in Scenario 2 are actually not better than those for Scenario 1. Possibly this is due to the rather unnatural data generation (using low degrees for the Laguerre polynomials only, leading to a sharp peak at the origin). Still, also Scenario 3 has a somewhat worse overall performance than Scenario 1; the model might be better at handling symmetric data.

\section{Real data application}\label{sec:realdata}
In this section, we apply our model to a real data example. We use data from the book of \cite{collett_modelling_2015}, that are discussed specifically in the context of dependent censoring there. These data have also been analysed in the dependent censoring setup of \cite{deresa_flexible_2020}. The description of the data that follows below, is based on Collett's Example 14.1 (pp. 461-463).

\subsection{Description of the data}
The data concern liver transplants in the UK. More specifically, the study cohort consists of 281 adult patients with primary biliary cirrhosis, a liver disease affecting mainly females, that were registered for receiving a transplant. However, because of a national shortage in livers available for transplant, some patients died while waiting for their transplant. The study was conducted in order to determine both the mortality rate on this waiting list as well as the impact of some factors on this rate; the UK Transplant Registry provided data on the time from first registration on the list (in the five-year period starting on 1 January 2006) to death. Some patients were removed from the list because a transplant was no longer possible due to their bad health condition. In their case, the removal from the list is considered as their death time. On the other hand, patients who did receive a transplant are considered censored at the time of the transplant, which happened for 217 out of 281 patients, inducing a high amount (about 77\%) of censoring. Since livers are allocated first to patients having the worst health condition, for those patients that are censored (receiving a liver), it is more likely that they would have had a shorter time to death than those that are uncensored. In other words, there is (positive) dependence in the censoring.

\begin{figure}[h]
	\begin{center}
	\subfigure[Time to death $T$ (+) or censoring $C$ ($\bullet$).]{
		\includegraphics[width=0.45\textwidth]{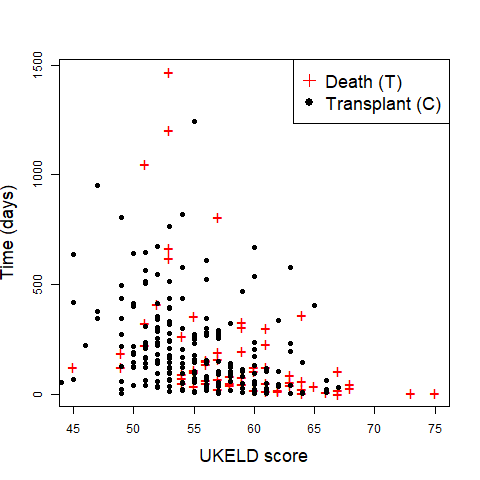}\label{fig:realdata}
	}
	\subfigure[Logarithmic time to death $\log(T)$.]{
		\includegraphics[width=0.45\textwidth]{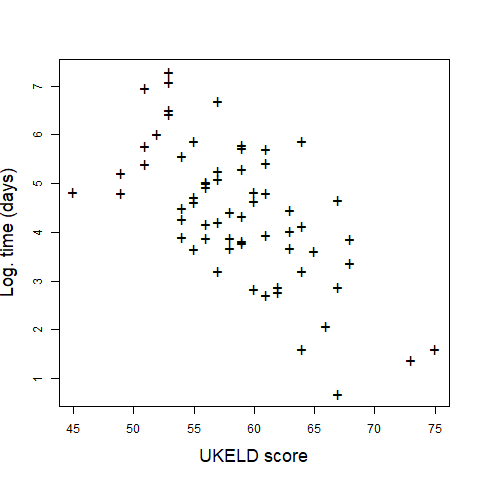}\label{fig:realdata:logT}
	}
	\end{center}
	\caption{Time to death (in days) while waiting for a liver transplant ($T$) or censoring due to receiving a transplant ($C$), versus the UKELD score. With a logarithmic scale for $T$, Figure (b) on the right shows a rather homoscedastic scenario and a linear trend can be observed.}
\end{figure}

The dataset contains information on some possible explanatory variables as well, namely the age (in years), sex (male = 1, female = 0), body mass index (BMI), and UK end-stage liver disease (UKELD) score. The latter is an indicator of the severity of the disease: the higher the UKELD score, the more severe the disease. As noted by \cite{deresa_flexible_2020}, it seems plausible that the time to death in this context depends on both the UKELD score and, conditional on a given UKELD score, on the time until receiving a transplant (see Figure \ref{fig:realdata}, adapted from their Figure 1).
Since both that paper as well as \cite{collett_modelling_2015} (Example 14.3, pp. 467-469) found in their analysis that only the UKELD score is significant, we omit the other explanatory variables in our analysis. In this way, the number of parameters to be estimated is reduced, such that we expect a more stable numerical performance compared to including superfluous explanatory variables. In order to enhance computational stability, our covariate $\tilde{X}$ is moreover a standardised version of the UKELD score. Recall that the full covariate $X = (1, \tilde{X})$ also includes a one for the intercept. Overall, this means that we fit the following model to the data:
\begin{equation}\label{eq:realdata:model}
	\begin{aligned}
		T &= \beta_0 + \beta_1\tilde{X} + \exp(\gamma_0 + \gamma_1\tilde{X}) \epsilon_T, \qquad &\epsilon_T \sim \eal, \\
		C &= \alpha_0 + \alpha_1\tilde{X} + \sigma_C \epsilon_C, &\qquad \epsilon_C \sim \mathcal{N}(0,1),
	\end{aligned}
\end{equation}
where $T$ denotes the logarithmic time to death (in days) while waiting for a liver transplant, and $C$ denotes (logarithmic) censoring due to receiving such a transplant. An overview of the logarithmic survival time $T$ versus the UKELD score is plotted in Figure \ref{fig:realdata:logT}; a linear trend can be observed, such that considering \eqref{eq:realdata:model} indeed makes sense. In the formulation of this model, the copula family as well as the heteroscedastic or homoscedastic nature is still open. These model components are varied over several models, the best of which is subsequently selected using AIC.

\subsection{Model results}
We fit several models to the data and then select the best one based on the AIC criterion. All copula families discussed above (i.e. Clayton, Gumbel and Frank) are considered, and so is the independence copula. Next, we also include an additional Frank copula model whose parameter $\theta$ is restricted such that the corresponding value of Kendall's tau lies between 0 and 1, called FrankPos. Both the nature of the data and the results obtained in \cite{deresa_flexible_2020} point in the direction of positive dependence between $T$ and $C$, so possibly better results are obtained when this extra information is supplied to the model. Figure \ref{fig:realdata:logT} suggests a rather homoscedastic behaviour. Nevertheless, we fit all models in the heteroscedastic setting as well for comparison. In the homoscedastic case, a fixed value for $\lambda$ between 0 and 1 is to be specified. We set this parameter to 0.5 in all homoscedastic cases initially. After selecting the optimal model according to AIC (round 1), only for this copula family we fitted some extra models (round 2) with different values for $\lambda$; recall from the simulations that in practice we do see performance differences, even if in theory any value for $\lambda$ should be equivalent.

\begin{table}[ht]
	\centering
	\caption{Fitted models (round 1) with the corresponding AIC scores. In the homoscedastic case, $\lambda$ equals $0.5$.}
	\label{tab:AIC:round1}
	\begin{tabular}{lccccc}
		\toprule
		& Indep & Gumbel & Clayton & Frank & FrankPos \\
		\midrule
		Hom. & 1157.20 & 1159.61 & 1158.18 &  1154.33 & 1155.09 \\
		Het. &  1157.47 & 1159.92 & 1156.04  & 1156.37 & 1158.10 \\ 
		\bottomrule
	\end{tabular}
\end{table}
Table \ref{tab:AIC:round1} above shows the resulting AIC values, both for the homoscedastic and the heteroscedastic case, in round 1. There doesn't seem to be much difference between the homoscedastic and the heteroscedastic models. Yet, the two best models according to AIC (Frank and FrankPos) are both on the homoscedastic side and, moreover, within the Frank family. Therefore, in round 2, we fitted the homoscedastic Frank and FrankPos models for some additional choices for $\lambda$. These results are displayed in Table \ref{tab:AIC:round2} below. Note that the ranking of all models in this table is the same on either side, regardless of the restriction on $\tau_K$, which is reassuring. Not so surprisingly, too extreme values for $\lambda$ (0.1 and 0.9) perform worst; Figure \ref{fig:realdata:logT} does not show any indication of the distribution for $T$ being this skewed. Within the more moderate range of values in the middle three columns, lower $\lambda$ lead to lower AIC values. Overall, we select the homoscedastic FrankPos model with $\lambda = 0.3$ as the best model according to AIC.

\begin{table}[ht]
	\centering
	\caption{Fitted models (round 2) with the corresponding AIC scores. Only in the second row, the copula parameter is restricted such that only positive dependence can be modelled.}
	\label{tab:AIC:round2}
	\begin{tabular}{lccccc}
		\toprule
		& $\lambda = 0.1$ & $\lambda = 0.3$ & $\lambda = 0.5$ & $\lambda = 0.7$ & $\lambda = 0.9$ \\
		\midrule
		Frank & 1158.21 & 1153.27 & 1154.33 & 1157.80 & 1159.73 \\
		FrankPos & 1159.86 & 1150.55 & 1155.09 & 1157.68 & 1162.12 \\
		\bottomrule
	\end{tabular}
\end{table}

\subsubsection{Parameter estimates}
Table \ref{tab:para:est} presents all parameter estimates in the selected model, as well as their counterpart in the best independence model, which is also the homoscedastic one. The variance term $\exp(\gamma_0 + \gamma_1 \tilde{X})$ in the selected homoscedastic case reduces to $\sigma_T = \exp(\gamma_0)$. A bootstrap standard error (BSE), estimated using 500 replications (in each of them determining the optimal Laguerre degrees for that particular replication), is also provided. Both results seem rather different; we quantify this using a likelihood ratio test (LRT). Use a superscript $d$ for the best dependent model and $i$ to denote the independent one. Write $q$ for the number of parameters in the model considered. Since the independence model does not have a dependence parameter but has, on the other hand, two extra Laguerre coefficients, it has one additional parameter compared to the FrankPos model. Computing the LRT statistic and comparing it to the critical chi-squared value corresponding to a significance level of $0.95$ and 1 degree of freedom (denoted $\chi^2_{1,95}$), we therefore get
\[\left(\text{AIC}_i - \text{AIC}_d\right) - 2 \left(q_i - q_d\right) = (1157.20 - 1150.55) - 2 = 5.65 > 3.84 \approx \chi^2_{1,95}.\]
Thus, we conclude that our FrankPos model with $\lambda = 0.3$ is indeed significantly better than the independence one. Also note the strong dependence ($\hat{\tau}_K \approx 0.61$), which is in line with the earlier findings of \cite{deresa_flexible_2020}.

\begin{table}[ht]
	\centering
	\caption{Parameter estimates (est.) and corresponding BSE based on 500 replications for the homoscedastic Frank model with $\lambda = 0.3$ and additional restriction for positive dependence, compared to those for the homoscedastic independence model with $\lambda = 0.5$. No BSE are included for the coefficients $\tilde{\phi}$ and $\phi$, as they are not present in every replication. \enquote*{--} means that the parameter is not present in the presented model (though, in the case of the Laguerre coefficients, they may be present in some of the bootstrap replications).}
	\label{tab:para:est}
	\begin{adjustbox}{width=\textwidth}
		\begin{tabular}{llccccccccccc}
			\toprule
			& & $\tau_K$ & $\beta_0$ & $\beta_1$ & $\sigma_T$ & $\tilde{\phi}_1$ & $\phi_1$ & $\phi_2$ & $\phi_3$ & $\alpha_0$ & $\alpha_1$ & $\sigma_C$ \\
			\midrule
			\multirow{2}{*}{FrankPos} & est. & 0.614 & 4.592 & -0.934 & 0.195 & -0.500 & -0.500 & -- & -- & 4.763 & -0.551 & 1.169 \\
			& BSE & 0.050 & 0.376 & 0.098 & 0.149 & & & & & 0.070 & 0.068 & 0.118 \\
			\midrule
			\multirow{2}{*}{Indep} & est. & -- & 5.963 & -1.057 & 0.288 & -0.140 & -0.309 & -2.423 & -0.568 & 4.942 & -0.374 & 1.255 \\
			& BSE & -- & 0.239 & 0.073 & 0.147 & & & & & 0.027 & 0.036 & 0.077 \\
			\bottomrule
		\end{tabular}
	\end{adjustbox}
\end{table}

\subsubsection{Quantile estimates}
The easiest way to interpret the parameter estimates from Table \ref{tab:para:est} is to consider quantile prediction for the (logarithmic) survival time. More precisely, on the level of the quantiles, model (\ref{eq:realdata:model}) in the homoscedastic case corresponds to
\[
Q_{T|X}(p|x) = \beta_0 + \beta_1\tilde{x} + \sigma_T Q_{\epsilon_T}(p)
\]
for the logarithmic survival $T$ (cf. \eqref{eq:qushift} \emph{supra}). For $p = 0.3$, the last term on the right vanishes (since also $\lambda$ equals 0.3) and we are left with $Q_{T|X}(0.3|x) = \beta_0 + \beta_1 \tilde{x}$ for the FrankPos model, where $\tilde{X}$ is the standardised UKELD score. Hence, for someone with the mean UKELD score ($\tilde{x} = 0$), the 0.3-quantile for the nonlogarithmic survival would be at $\exp(\beta_0) \approx \exp(4.59) \approx 99$ days. On the other hand, for someone with a UKELD score that is one standard deviation to the right (approximately 5.79 UKELD units), the covariate $\tilde{x}$ would be 1, such that the 0.3 quantile is now only at $\exp(\beta_0 + \beta_1) \approx \exp(4.59 - 0.93) \approx 39$ days. The UKELD score thus drastically reduces the predicted quantile of the survival time. For comparison, we compute the corresponding 0.3-quantiles for the independence model. For the mean UKELD score ($\tilde{x} = 0$), we get approximately $\exp(5.56) \approx 261$ days, and for $\tilde{x} = 1$, the 0.3-quantile estimate under independence is at about $\exp(4.51) \approx$ 91 days. The relative decrease for 1 standard deviation in UKELD units is roughly speaking similar (in both cases a factor between 2.5 and 3), but the values themselves are far from being the same. Similar computations can be made for other quantiles.

\begin{table}[h]
	\centering
	\caption{Quantile estimates for the best FrankPos model (i.e. homoscedastic Frank with $\lambda = 0.3$ and extra restriction for positive dependence), compared to those for the homoscedastic independence model; the estimated bootstrap standard error based on 500 replications is included in parentheses.}
	\label{tab:qu:est}
	\begin{tabular}{ll@{\hskip 25pt}c@{\hskip 25pt}c@{\hskip 25pt}c}
		\toprule
		& & \multicolumn{3}{c}{UKELD score} \\
		& & $50$ & $60$ & $70$ \\
		\midrule
		\multirow{3}{*}{FrankPos}
		& $p = 0.3$ & 5.674 (0.461) & 3.806 (0.325) & 1.939 (0.275) \\
		& $p = 0.5$ & 6.373 (0.500) & 4.505 (0.370) & 2.638 (0.316) \\
		& $p = 0.7$ & 7.063 (0.611) & 5.195 (0.502) & 3.328 (0.455) \\
		\midrule
		\multirow{3}{*}{Indep}
		& $p = 0.3$ & 6.790 (0.187) & 4.675 (0.102) & 2.561 (0.170) \\
		& $p = 0.5$ & 7.188 (0.278) & 5.074 (0.228) & 2.959 (0.264) \\
		& $p = 0.7$ & 8.087 (0.392) & 5.973 (0.394) & 3.858 (0.448) \\
		\bottomrule
	\end{tabular}
\end{table}

Supplement S5 contains an overview of the estimated quantiles for $p$ ranging from 0.1 to 0.9 and some UKELD scores, with a brief discussion. The trend is very clear in the summary Table \nolinebreak \ref{tab:qu:est} included above: the independence model has a larger estimate for each quantile. Thus, ignoring the dependence in the data leads to overestimating the time to death. The differences between both models are moreover quite large. Recall that all quantiles reported in Table \ref{tab:qu:est} are for the logarithmic survival time. Therefore, for the survival time itself, these (already not so small) differences get huge very rapidly. The median of the survival time for a UKELD score of $60$, for instance, is approximately $160$ days in the independent censoring model, compared to only $90$ under our dependent model -- this means that the median according to the FrankPos model is only about 56\% of what the independence model predicts. The difference between the two models gets even larger for more extreme quantiles and UKELD scores. We thus conclude that the difference in the model parameters translates into quite different quantile predictions as well; the dependent nature of the censoring cannot be ignored in this data application.

\section{Concluding remarks} \label{sec:concl}
With this article, we contribute to the very limited literature on quantile regression under dependent censoring. Previous related research was often restricted to a semicompeting risks context (e.g. \cite{li_quantile_2015} and references therein) or, within a competing risks framework, focused on cause-specific cumulative incidence functions, as in \cite{peng_nonparametric_2007, peng_competing_2009}. The work on dependent censoring doing inference for the \emph{net} quantity $T|X$ can be seen as an interesting alternative. In our work, the partial identification of \cite{fan_partial_2018} is upgraded to full parameter identification. On the other hand, the modelling assumptions in \cite{ji_analysis_2014} of a fixed copula and (transformed) linear quantiles for both margins, are relaxed to a flexible association, and linear quantiles for $T$ only, moreover not necessarily on the quantile level of interest (cf. Section \ref{sec:model:multiqu}).

Such relaxation possibilities are due to the distributional flavour in our quantile regression, that enables application of \cite{czado_dependent_2023}. Given the existing extensions of that article to models where only one of the marginal models is fully parametric \cite{deresa_copula_2024}, also for our work a possible extension could consist in modelling one of the margins non- or semiparametrically. If one would like to preserve the current quantile focus introduced by using the specific EAL distribution with its Laplace basis and Laguerre enrichment (introduced in \cite{kreiss_VK_submitted}), the non- or semiparametric component would have to be the censoring, or a more flexible variant of the currently linear term $X^\top \beta$. Identifiability, however, would not be obvious, as the proof relies on easy identification of $C$ to determine $T$'s parameters only afterwards; possibly in this setup we cannot get much closer to the identifiability borders posed by \cite{tsiatis_nonidentifiability_1975}. More straightforward generalisations could be the inclusion of administrative censoring and a covariate-dependent copula parameter to accomodate for inhomogeneous associations, for which the proof seems to be readily generalisable. Yet, it is only fair to mention that all these generalisations would further increase the number of parameters and entail an even more heavy numerical estimation procedure. Perhaps, a more plausible alternative path would consist in combining the somewhat simpler two-piece distributions of \cite{ewnetu_flexible_2023, ewnetu_two-piece_2023}, that share our quantile orientation, with the parametric copula approach of \cite{czado_dependent_2023} and a more flexible marginal model for the censoring.

We finally point out that the distributional touch in our quantile regression takes away identifiability issues for extreme quantiles (cf. Remark \ref{rem:model}\ref{rem:model:tautrunc}), and, moreover, leads to the shift concept of Section \ref{sec:model:multiqu}. Its applicability is not merely restricted to EAL distributions; the families of two-piece distributions in \cite{gijbels_quantile-based_2019, gijbels_semiparametric_2021} and \cite{ewnetu_two-piece_2023}, for instance, also allow for this shift approach, and they already noted the corresponding preclusion of crossing quantile curves. However, none of them mentions the large computational advantage of having to perform maximum likelihood optimisation only once, nor the nice theoretical property of more flexible linearity assumptions. Further exploring the applicability and advantages of both this shift approach and such less standard \enquote*{distributional quantile regression} in other contexts might be an interesting topic for future research.


\clearpage
\appendix

\renewcommand{\thefigure}{A.\arabic{figure}}
\setcounter{figure}{0}
\renewcommand{\thetable}{A.\arabic{table}}
\setcounter{table}{0}
\renewcommand{\thesection}{A\arabic{section}}
\renewcommand{\thesubsection}{\thesection.\arabic{subsection}}
\renewcommand{\thesubsubsection}{\thesubsection.\arabic{subsubsection}}

\section{Appendix: Assumption verification: auxiliary lemmata}\label{app:idf:auxlem}
In this first appendix we provide auxiliary lemmata that can be used to verify the validity of the identifiability assumptions in Section \ref{sec:idf:ass}. All proofs can be found in Supplement S1.2; the proofs of Lemma \ref{lem:hfun:lim} and Lemma \ref{lem:hfun:idf:theta:clayton} are inspired by the proof of Theorem 3 and Theorem 4, respectively, in \cite{czado_dependent_2023}.

To start, whether or not assumption \ref{ass:idf:LC} and \ref{ass:idf:LT} are satisfied depends both on the limiting behaviour of the $h$-functions involved, and on the specific distribution for $C$. To assess the validity of assumption \ref{ass:idf:LC:hfun} and \ref{ass:idf:LT:hfun}, depending on the copula family to which $\cpl_\theta(\cdot, \cdot)$ belongs, the following lemma comes in handy.

\begin{Lemma}\label{lem:hfun:lim}
	Let $x \in R$ and $(\theta, \theta_T, \theta_C) \in \Theta \times \Theta_T \times \Theta_C$ be any set of parameters. Then the equations
	\begin{subequations}
		\begin{eqnarray}
			\lim\limits_{y \to L_C} h_{T|C; \theta}(F_{T|X; \theta_T}(y|x)| F_{C|X; \theta_C}(y|x)) &=& 0 \label{eq:lem:hfun:lima} \\
			\lim\limits_{y \to - \infty} h_{C|T; \theta}(F_{C|X; \theta_C}(y|x)| F_{T|X; \theta_T}(y|x)) &=& 0 \label{eq:lem:hfun:limb}
		\end{eqnarray}
	\end{subequations}
	corresponding to the copula function $\cpl_\theta(\cdot, \cdot)$ hold true in the following situations:
	\begin{enumerate}[(i)]
		\item Independence copula: \eqref{eq:lem:hfun:lima} holds for $L_C = - \infty$; also \eqref{eq:lem:hfun:limb} holds.
		\item Frank copula: \eqref{eq:lem:hfun:lima} holds for $L_C = - \infty$; also \eqref{eq:lem:hfun:limb} holds.
		\item Gumbel copula: whenever $\theta > 1$ (the case $\theta = 1$ corresponds to the independence copula),
		\begin{itemize}
			\item \eqref{eq:lem:hfun:lima} holds for $L_C = M_r$ (the upper bound of $C|X$'s support) if $M_r < \infty$, and
			\item \eqref{eq:lem:hfun:limb} only holds provided that
			\begin{equation}\label{ass:idf:G3}\tag{G3}
				\lim\limits_{y \to - \infty} \frac{\log F_{C|X; \theta_C}(y|x)}{\log F_{T|X; \theta_T}(y|x)} \in (0, \infty].
			\end{equation}
		\end{itemize}
		\item Clayton copula: whenever condition \ref{ass:idf:B1} and \ref{ass:idf:B2} are satisfied, 
		equation \eqref{eq:lem:hfun:lima} holds for $L_C = - \infty$, but \eqref{eq:lem:hfun:limb} does \emph{not} hold.
	\end{enumerate}
	
\end{Lemma}
\noindent
Once these conditions are verified, one can check whether also \ref{ass:idf:LC:asyidf} holds for this value of $L_C$.
\bigbreak
Also the validity of assumption \ref{ass:idf:coppara}, crucial for the identification of the copula parameter $\theta$, is not obvious in most cases. Once again, for some of the common Archimedean copula families, a lemma can be used to verify that \ref{ass:idf:coppara} holds true.

\begin{Lemma}\label{lem:hfun:idf:theta}
	Consider two candidate pairs $(\theta_1, \theta_{T_1}), (\theta_2, \theta_{T_2}) \in \Theta \times \Theta_T$ and let $\theta_C$ be any value in $\Theta_C$. Suppose that for the $h_{T|C}$-functions corresponding to the copula functions $\cpl_{\theta_1}(\cdot, \cdot)$ and $\cpl_{\theta_2}(\cdot, \cdot)$ belonging to the same Archimedean copula family,
	\begin{equation}\label{eq:lem:hfun}
		1 = \frac{h_{T|C; \theta_1}(F_{T|X; \theta_{T_1}}(y|x) | F_{C|X; \theta_{C}}(y|x))}{h_{T|C; \theta_2}(F_{T|X; \theta_{T_2}}(y|x) | F_{C|X; \theta_{C}}(y|x))}
	\end{equation}
	for any $x$ and corresponding $y \in U_C(x)$. Then $\theta_1 = \theta_2$ is implied whenever:
	\begin{enumerate}[(i)]
		\item $\cpl_\theta(\cdot, \cdot)$ is the independence copula.
		\item $\cpl_\theta(\cdot, \cdot)$ belongs to the Frank family and in addition, assumption \ref{ass:idf:F1} holds:
		\begin{enumerate}[label = (F\arabic{enumii})]
			\item \label{ass:idf:F1} For at least one $x \in R$, $U_C(x)$ contains an interval $I_x = (-\infty, A_x)$ for some $A_x \in \R$.
		\end{enumerate}
		\item $\cpl_\theta(\cdot, \cdot)$ belongs to the Gumbel family and in addition, assumption \ref{ass:idf:GumbInt} holds:
		\begin{enumerate}[label = (G\arabic{enumii})]
			\setcounter{enumii}{3}
			\item \label{ass:idf:GumbInt} For at least one $x \in R$, $U_C(x)$ contains an interval $I_x = (A_x, M_r)$ for some $A_x \in \R$, where the endpoint $M_r$ of the support of $C|X$ is moreover finite.
		\end{enumerate}
	\end{enumerate}
\end{Lemma}

The Clayton case has been excluded from the previous lemma, since a stronger result can be formulated in this case. This is useful: as we have seen in Lemma \ref{lem:hfun:lim}, we can never have both assumption \ref{ass:idf:LC} and assumption \ref{ass:idf:LT} satisfied, that are key assumptions in the proof to identify $\theta_C$ and $\theta_T$, respectively. On the other hand, the proof method of Theorem 4 in \cite{czado_dependent_2023} shows that equality of the $h_{T|C}$-functions as in the statement of Lemma \ref{lem:hfun:idf:theta} above can be used to also obtain some information for $T$ already. In this way, even using only one type of $h$-function yields a sufficient amount of information to perform all necessary steps in the identifiability proof.

\begin{Lemma}\label{lem:hfun:idf:theta:clayton}
	Consider two candidate pairs $(\theta_1, \theta_{T_1}), (\theta_2, \theta_{T_2}) \in \Theta \times \Theta_T$ and let $\theta_C$ be any value in $\Theta_C$. Suppose that for the $h_{T|C}$-functions corresponding to the Clayton copula functions $\cpl_{\theta_1}(\cdot, \cdot)$ and $\cpl_{\theta_2}(\cdot, \cdot)$, equation (\ref{eq:lem:hfun}) holds for any $x$ and corresponding $y \in U_C(x)$. Suppose moreover that assumption \ref{ass:idf:B3} holds. Then $\theta_1 = \theta_2$, and moreover for any $x \in R$:
	\begin{equation*}
		1 = \lim\limits_{y \to - \infty}\frac{f_{T|X; \theta_{T_1}}(y|x)}{f_{T|X; \theta_{T_2}}(y|x)}.
	\end{equation*}
\end{Lemma}

As is apparent from the identification proof in Appendix \ref{app:proof:idfthm}, this allows replacing assumptions \ref{ass:idf:LC}-\ref{ass:idf:coppara} by the alternative assumptions \ref{ass:idf:B1}-\ref{ass:idf:B4}.
\bigbreak
Finally, we state the very simple, general Lemma \ref{lem:hTC:inj} that is useful in one of the steps of the identifiability proof; its validity for the main Archimedean copulas under consideration is given in the subsequent Lemma \ref{lem:hTC:inj:satisfied}.

\begin{Lemma}\label{lem:hTC:inj}
	Suppose $\cpl_\theta(\cdot, \cdot)$ is an Archimedean copula with generator $\psi_{\theta}(t)$ for $t \in [0,1]$. Assume that on the unit interval, $\psi^{-1}(\cdot)$ exists and $\psi_\theta'(\cdot)$ is an injective function. Then for any $v \in [0,1]$, the function $h_{T|C; \theta}(\cdot | v)$ is injective on $[0,1]$.
\end{Lemma}

\begin{Lemma}\label{lem:hTC:inj:satisfied}
	The conditions of Lemma \ref{lem:hTC:inj} are satisfied for the following Archimedean copula families, and for the associated domains for the copula parameter $\theta$:
	\begin{enumerate}[(i)]
		\item the independence copula
		\item the Frank copula with $\theta \in \R \setminus \{0\}$
		\item the Gumbel copula with $\theta \geq 1$
		\item the Clayton copula with $\theta \in (-1, \infty) \setminus \{0\}$, but not for $\theta = -1$.
	\end{enumerate}
\end{Lemma}

\section{Appendix: Identifiability proof (Theorem \ref{thm:idf} and \ref{thm:idf:clayton})}\label{app:proof:idfthm}
This appendix contains the proof of the identifiability statement formulated in Theorem \ref{thm:idf} above. All its steps are valid regardless of the copula family and marginal distributions under consideration, whenever assumptions \ref{ass:idf:posdef}-\ref{ass:idf:S2} hold. Assumption \ref{ass:idf:posdef} does not explicitly appear in this proof, but is used to prove the identifiability of $\beta$ depending on the status of $\lambda$, cf. Supplement S1.1, and this case distinction is assumed in the proof. Verifying the validity of assumptions \ref{ass:idf:LC}-\ref{ass:idf:coppara} is highly copula-dependent; the necessary auxiliary lemmata are proved in Supplement S1.2. The idea of using limiting information of the $h$-functions is taken from the proof of Theorem 1 in \cite{czado_dependent_2023}. However, our proof requires a way lengthier reasoning since $\theta_T$ in our case cannot be identified from its asymptotic behaviour only, contrary to their assumptions (i.e. the hypothetical analogue of assumption \ref{ass:idf:LC:asyidf} for $T$, cannot be true for our EAL distributed $T$).
\bigbreak
The proof presented here uses the validity of assumptions \ref{ass:idf:posdef}-\ref{ass:idf:S2}. In the case of the Clayton copula, we work with the alternative set of assumptions \ref{ass:idf:B1}-\ref{ass:idf:B4} and \ref{ass:idf:posdef}-\ref{ass:idf:gamidf}, \ref{ass:idf:S0}-\ref{ass:idf:S2} instead, once again inspired by \cite{czado_dependent_2023}, now Theorem 4. The proof being highly similar, we just point out the few differences in the final subsection.

\subsection{Identification of $C$'s parameters $\theta_C$}\label{app:proof:idfthm:C}
Recall that the likelihood contribution for a censored observation is of the form (\ref{eq:lik:contr}) with $\delta = 0$, such that for any $x$ and $y$ with positive censoring probability, there must hold
\begin{equation}\label{eq:proof:gumbel:cens:llh}
	1 = \frac{f_{C|X; \theta_{C_1}}(y|x)}{f_{C|X; \theta_{C_2}}(y|x)} \cdot \frac{1 - h_{T|C; \theta_1}(F_{T|X; \theta_{T_1}}(y|x) | F_{C|X; \theta_{C_1}}(y|x))}{1 - h_{T|C; \theta_2}(F_{T|X; \theta_{T_2}}(y|x) | F_{C|X; \theta_{C_2}}(y|x))}.
\end{equation}
Now take any $x \in R$ and a corresponding set of values $y \to L_C$ as in assumption \ref{ass:idf:LC:seq}, such that the above equation holds for any such pair $(y,x)$. In view of assumption \ref{ass:idf:LC:hfun}, we thus get
\begin{equation*}
	\lim\limits_{y \to L_C} \frac{f_{C|X; \theta_{C_1}}(y|x)}{f_{C|X; \theta_{C_2}}(y|x)} = 1
\end{equation*}
for any $x \in R$, which, by assumption \ref{ass:idf:LC:asyidf}, identifies $\theta_C$.

\subsection{Identification of the copula parameter $\theta$}\label{app:proof:idfthm:cop}
As $f_{C|X}(\cdot| \cdot)$ is now fully identified, (\ref{eq:proof:gumbel:cens:llh}) implies the following equality for any tuple $(y,x)$ belonging to $U_C(x) \times \{x\}$ for any $x$ (hence with positive censoring probability):
\begin{equation}\label{eq:proof:gumbel:hbasis}
	1 = \frac{h_{T|C; \theta_1}(F_{T|X; \theta_{T_1}}(y|x) | F_{C|X; \theta_{C}}(y|x))}{h_{T|C; \theta_2}(F_{T|X; \theta_{T_2}}(y|x) | F_{C|X; \theta_{C}}(y|x))}.
\end{equation}
The fact that this implies equality of the copula parameter candidates $\theta_1$ and $\theta_2$ is now precisely the content of assumption \ref{ass:idf:coppara}.

\subsection{Identification of $T$'s parameters $\theta_T$}\label{app:proof:idfthm:T}
The likelihood contribution for an uncensored observation being of the form given in (\ref{eq:lik:contr}) with $\delta = 1$, for any $x$ and $y$ with positive probability of being observed (i.e. $y \in U_T(x)$), one has
\begin{equation*}
	1 = \frac{f_{T|X; \theta_{T_1}}(y|x)}{f_{T|X; \theta_{T_2}}(y|x)} \cdot \frac{1 - h_{C|T; \theta_1}(F_{C|X; \theta_{C_1}}(y|x) | F_{T|X; \theta_{T_1}}(y|x))}{1 - h_{C|T; \theta_2}(F_{C|X; \theta_{C_2}}(y|x) | F_{T|X; \theta_{T_2}}(y|x))}
\end{equation*}
(the parameters $\theta$ and $\theta_C$ are, actually, identified already, but this information is not necessary here). Using assumption \ref{ass:idf:LT:seq}, take any $x \in R$ and a corresponding sequence $y \to - \infty$ within $U_T(x)$, then the previous equality together with the limiting behaviour of the $h$-function in assumption \ref{ass:idf:LT:hfun} yields for any $x \in R$:
\begin{equation}\label{eq:proof:gumbel:tlim}
	\lim\limits_{y \to - \infty} \frac{f_{T|X; \theta_{T_1}}(y|x)}{f_{T|X; \theta_{T_2}}(y|x)} = 1.
\end{equation}
Combining \eqref{eq:Tdistr} with the density function of the EAL distribution (\ref{eq:dens:eal}) and the limit in (\ref{eq:proof:gumbel:tlim}) above, we get that for any value $X = x$, and using $\sigma_i(x)$ as shorthand notation for $\sigma(x; \gamma_i)$:
\begin{multline*}
	1 = \frac{\lambda_1 (1 - \lambda_1)}{\lambda_2 (1 - \lambda_2)}
	\cdot \left(\frac{\|\tilde{\phi}_2\|}{\|\tilde{\phi}_1\|}\right)^2
	\cdot \frac{\sigma_2(x)}{\sigma_1(x)}
	\cdot \lim\limits_{y \to -\infty} \frac{\left(\sum\limits_{k = 0}^{\tilde{m}} \tilde{\phi}_{1k} L_k\left(\frac{\lambda_1 - 1}{\sigma_1(x)}(y - x^\top \beta_1)\right)\right)^2}{\left(\sum\limits_{k = 0}^{\tilde{m}} \tilde{\phi}_{2k} L_k\left(\frac{\lambda_2 - 1}{\sigma_2(x)}(y - x^\top \beta_2)\right)\right)^2}
	\cdot \\
	\lim\limits_{y \to -\infty} \exp\left(\left(\frac{\lambda_2 - 1}{\sigma_2(x)} - \frac{\lambda_1 - 1}{\sigma_1(x)}\right) y + x^\top\left(\frac{(\lambda_1 - 1)\beta_1}{\sigma_1(x)} - \frac{(\lambda_2 - 1)\beta_2}{\sigma_2(x)}\right)\right).
\end{multline*}
The limit on the first line is determined by its highest polynomial degrees, so
\begin{multline}\label{eq:proof:gumbel:split}
	1 = \frac{\lambda_1}{\lambda_2}
	\cdot \left(\frac{\|\tilde{\phi}_2\|}{\|\tilde{\phi}_1\|}\right)^2
	\cdot \left(\frac{\tilde{\phi}_{1m}}{\tilde{\phi}_{2m}}\right)^2
	\cdot \left(\frac{\frac{\lambda_1 - 1}{\sigma_1(x)}}{\frac{\lambda_2 - 1}{\sigma_2(x)}}\right)^{2\tilde{m} + 1}
	\cdot \exp \left(x^\top\left(\frac{(\lambda_1 - 1)\beta_1}{\sigma_1(x)} - \frac{(\lambda_2 - 1)\beta_2}{\sigma_2(x)}\right)\right)
	\cdot \\
	\lim\limits_{y \to -\infty} \exp\left(\left(\frac{\lambda_2 - 1}{\sigma_2(x)} - \frac{\lambda_1 - 1}{\sigma_1(x)}\right) y\right).
\end{multline}
This can only happen if the remaining limit is nonzero and finite, i.e. if
\begin{equation}\label{eq:proof:gumbel:lamsig}
	\frac{\lambda_1 - 1}{\sigma_1(x)} = \frac{\lambda_2 - 1}{\sigma_2(x)},
\end{equation}
and, moreover, since (\ref{eq:proof:gumbel:split}) must hold for arbitrary $x$, the argument of the exponential function must be constant. After substitution of (\ref{eq:proof:gumbel:lamsig}), this in particular implies that there exists some constant $K$ for which
\begin{equation*}
	K = (\beta_1 - \beta_2)^\top \frac{X}{\sigma_1(X)}.
\end{equation*}

\subsubsection{Identification of $\beta, \lambda$ and $\gamma$}
As explained in Supplement S1.1, whether the previous result allows identification of $\beta$ depends on $\sigma(X; \gamma)$. Therefore, case distinction is made for the case where $\beta$ can be identified (and hence $\lambda$ is kept variable) versus the case where it cannot be, without extra information (hence $\lambda$ is fixed), see also Section \ref{sec:model:multiqu:parastatus:lam}. In the first case, we now identify $\lambda$ based on our identified $\beta$, while in the second case, we identify $\beta$ based on our information for $\lambda$. In either case, the following identity is useful:
\begin{equation}\label{eq:proof:gumbel:FTonWC}
	\forall x, \forall y \in U_C(x): \quad F_{T|X; \theta_{T_1}}(y|x) = F_{T|X; \theta_{T_2}}(y|x),
\end{equation}
which, having identified the parameters of both $C$ and the copula already, we can deduce from the equality of the $h$-functions (\ref{eq:proof:gumbel:hbasis}) for the corresponding tuples with positive censoring probability, using Lemma \ref{lem:hTC:inj}.
\bigbreak
\noindent
\textbf{Case (i): $\lambda$ variable}.
\bigbreak
By assumption, this means that $\beta$ could already be identified; say it equals $\beta^0$. Assumption \ref{ass:idf:S1} then in particular implies that the set $S(\beta^0)$ is nonempty; taking any such $x \in S(\beta^0)$, assumption \ref{ass:idf:S0} implies that $x^\top \beta^0 \in U_C(x)$ and hence (\ref{eq:proof:gumbel:FTonWC}) can be applied to obtain
\begin{equation*}
	F_{T|X; \theta_{T_1}}(x^\top \beta^0|x) = F_{T|X; \theta_{T_2}}(x^\top \beta^0|x).
\end{equation*}
By the model construction, $x^\top \beta^0$ is the $\lambda$-th quantile, hence the left-hand side equals $\lambda_1$ while on the right we get $\lambda_2$. Thus $\lambda_1 = \lambda_2$.

\bigbreak
\noindent
\textbf{Case (ii): $\lambda$ fixed}.
\bigbreak
By assumption \ref{ass:idf:S2}, for any point $x$ in the intersection $S(\beta_1) \cap S(\beta_2)$, we can once again apply (\ref{eq:proof:gumbel:FTonWC}), to get
\begin{equation*}
	F_{T|X; \theta_{T_1}}(x^\top \beta_1|x) = F_{T|X; \theta_{T_2}}(x^\top \beta_1|x)
\end{equation*}
(note that $y = x^\top \beta_1$ twice, while the indices for $\theta_T$ differ). On the other hand, now using the opposite reasoning from before (starting from the fact that the $\lambda$-th quantile of $F_{T|X; \theta_i}$ equals $x^\top\beta_i$ by construction),
\begin{equation*}
	F_{T|X; \theta_{T_1}}(x^\top \beta_1|x) = \lambda = F_{T|X; \theta_{T_2}}(x^\top \beta_2|x).
\end{equation*}
Combining these last two equalities yields $F_{T|X; \theta_{T_2}}(x^\top \beta_1|x) = F_{T|X; \theta_{T_2}}(x^\top \beta_2|x)$. Since $T$ has a strictly increasing distribution function, injectivity implies that $x^\top(\beta_1 - \beta_2) = 0$ for all such $x$. Assumption \ref{ass:idf:S2} now guarantees that the set of points $x$ for which the above reasoning holds, is nonempty and contains enough information for $\beta$ to be identified. Indeed,
\begin{equation*}
	\beta_{10} - \beta_{20} = \tilde{X}^\top(\tilde{\beta}_2 - \tilde{\beta}_1)
\end{equation*}
on $S(\beta_1) \cap S(\beta_2)$, hence the right-hand side is constant and the positive definiteness assumption implies $\tilde{\beta}_1 = \tilde{\beta}_2$ which, in turn, implies that the expression on the left of the above equation equals zero as well; $\beta_1 = \beta_2$ follows.
\bigbreak
In either case, $\lambda$ is now identified. Equation (\ref{eq:proof:gumbel:lamsig}) then identifies $\sigma_i(X) = \sigma(X; \gamma_i)$ and in turn, by assumption \ref{ass:idf:gamidf}, its parameter $\gamma$.

\subsubsection{Identification of $\tilde{\phi}$ and $\phi$}
Finally, to identify the last parameters of $T$, we use the fact that two polynomials can be equal on an entire interval only if their coefficients are all identical. Using the value $x^0$ and its corresponding interval $I$ from assumption \ref{ass:idf:S1}, equation (\ref{eq:proof:gumbel:FTonWC}) implies that $F_{T|X; \theta_1}(y|x^0) = F_{T|X; \theta_2}(y|x^0)$ for all $y \in I$. Equality on the entire interval enables differentiation, such that the same equality now holds for the corresponding density functions on $I$. As $(x^0)^\top \beta^0$ (with $\beta^0$ the identified value for $\beta$) is an internal point, $I$ contains both a subinterval $I_- \subset I$ lying below $(x^0)^\top \beta^0$ and an $I_+ \subset I$ above this point. These intervals can be used to identify $\tilde{\phi}$ and $\phi$, respectively. The reasoning for $I_-$ and $\tilde{\phi}$ is outlined below; the one for $I_+$ and $\phi$ is analogous.
\bigbreak
Since all other parameters of $T$ are identified, equality of the density functions on some interval $I_-$ below $(x^0)^\top \beta^0$ is equivalent to the existence of some subinterval $J_- \subset \R^-$ of the negative real line for which
\begin{eqnarray*}
	f_{\eal, 1}(u) &=& f_{\eal, 2}(u), \qquad \forall u \in J_- \\
	\Iff \left(\sum_{k = 0}^{\tilde{m}} \frac{\tilde{\phi}_{1k}}{ \|\tilde{\phi_1}\|} L_k ((\lambda-1) u)\right)^2 &=& \left(\sum_{k = 0}^{\tilde{m}} \frac{\tilde{\phi}_{2k}}{ \|\tilde{\phi_2}\|} L_k ((\lambda-1) u)\right)^2, \qquad \forall u \in J_-.
\end{eqnarray*}
As equality of polynomials on an interval suffices to identify the coefficients, this implies that
\[\frac{\tilde{\phi}_{1k}}{ \|\tilde{\phi_1}\|} = \frac{\tilde{\phi}_{2k}}{ \|\tilde{\phi_2}\|}, \qquad \forall k.\]
(Note that we start from a linear combination of polynomials, so for each power of $u$, we get linear combinations of the polynomial coefficients of the Laguerre polynomials, with the $\frac{\tilde{\phi}_k}{\|\tilde{\phi}\|}$ as weights. It is not immediately clear that identifiability of these combinations also implies identifiability of the distinct $\tilde{\phi}$. However, when imposing $\tilde{\phi}_{10} = 1 = \tilde{\phi}_{20}$, a (backwards) recursive reasoning can be applied to show that the above equality indeed identifies the $\tilde{\phi}$.)
\qed

\subsection{Adaptations for the Clayton case (Theorem \ref{thm:idf:clayton})}
The reasoning in Appendix \ref{app:proof:idfthm:C} for the identification of $C$'s parameters remains valid, since \ref{ass:idf:B3} and \ref{ass:idf:B4} correspond to \ref{ass:idf:LC:seq} and \ref{ass:idf:LC:asyidf}, respectively (for $L_C = - \infty$), and \ref{ass:idf:LC:hfun} has been replaced by the result of Lemma \ref{lem:hfun:lim}, satisfied under \ref{ass:idf:B1}-\ref{ass:idf:B2}. Next, Lemma \ref{lem:hfun:idf:theta:clayton}, using assumption \ref{ass:idf:B3}, can be used to replace both the identification of the copula parameter in Appendix \ref{app:proof:idfthm:cop}, and the first part of Appendix \ref{app:proof:idfthm:T}, to arrive at the limiting behaviour (\ref{eq:proof:gumbel:tlim}). The remaining part of the identification of $T$'s parameters makes use of \ref{ass:idf:S0}-\ref{ass:idf:S2} and \ref{ass:idf:posdef}-\ref{ass:idf:gamidf} only.
\qed

\end{document}